\documentclass[UTF-8,reqno]{amsart}
\usepackage{enumerate}
\usepackage{amssymb,url,color, booktabs}
\usepackage{mathrsfs}

\usepackage[nobysame]{amsrefs}
\BibSpec{article}{%
+{}{\PrintAuthors} {author}
+{,}{ \textrm} {title}
+{.}{ \textit} {journal}
+{,}{ \textbf} {volume}
+{}{ \parenthesize} {date}
+{,}{ } {pages}
+{.}{ arXiv:} {eprint}
+{.}{} {transition}
}
\BibSpec{book}{%
+{}{\PrintAuthors} {author}
+{,}{ \textit} {title}
+{.}{ \textrm} {series} 
+{,}{ Vol.} {volume} 
+{.}{ } {publisher}
+{,}{ } {date}
+{.}{} {transition}
}
\usepackage{color}
\usepackage[colorlinks=true]{hyperref}
\hypersetup{
    linkcolor=blue,          
    citecolor=red,        
    filecolor=blue,      
    urlcolor=cyan  
}

\numberwithin{equation}{section}

\newcommand{\be}{\begin{eqnarray}}
\newcommand{\ee}{\end{eqnarray}}
\newcommand{\ce}{\begin{eqnarray*}}
\newcommand{\de}{\end{eqnarray*}}
\newtheorem{theorem}{Theorem}[section]
\newtheorem{lemma}[theorem]{Lemma}
\newtheorem{remark}[theorem]{Remark}
\newtheorem{definition}[theorem]{Definition}
\newtheorem{proposition}[theorem]{Proposition}
\newtheorem{Examples}[theorem]{Example}
\newtheorem{corollary}[theorem]{Corollary}

\def\eps{\varepsilon}

\def\e{\mathrm{e}}

\def\p{\partial}

\def\[{{\Big[}}
\def\]{{\Big]}}
\def\<{{\langle}}
\def\>{{\rangle}}
\def\({{\Big(}}
\def\){{\Big)}}

\def\bx{{\mathbf{x}}}
\def\tr{\mathrm {tr}}

\def\dif{{\mathord{{\rm d}}}}

\def\min{{\mathord{{\rm min}}}}

\def\={&\!\!=\!\!&}

\def\cB{{\mathcal B}}

\def\cF{{\mathcal F}}

\def\cL{{\mathcal L}}
\def\cM{{\mathcal M}}
\def\cN{{\mathcal N}}

\def\cP{{\mathcal P}}

\def\mC{{\mathbb C}}
\def\mD{{\mathbb D}}
\def\mE{{\mathbb E}}

\def\mI{{\mathbb I}}

\def\mN{{\mathbb N}}

\def\mP{{\mathbb P}}
\def\mQ{{\mathbb Q}}
\def\mR{{\mathbb R}}
\def\mS{{\mathbb S}}

\def\mX{{\mathbb X}}

\def\bP{{\mathbf P}}

\def\1{{\mathbf{1}}}

\def\sA{{\mathscr A}}
\def\sB{{\mathscr B}}

\def\sI{{\mathscr I}}

\def\sL{{\mathscr L}}
\def\sM{{\mathscr M}}
\def\sN{{\mathscr N}}

\def\sQ{{\mathscr Q}}

\def\geq{\geqslant}
\def\leq{\leqslant}

\def\div{\mathord{{\rm div}}}

\def\eps{\varepsilon}

\def\e{\mathrm{e}}

\def\p{\partial}

\def\[{{\Big[}}
\def\]{{\Big]}}
\def\<{{\langle}}
\def\>{{\rangle}}
\def\({{\Big(}}
\def\){{\Big)}}

\def\bx{{\mathbf{x}}}
\def\tr{\mathrm {tr}}

\def\dif{{\mathord{{\rm d}}}}

\def\min{{\mathord{{\rm min}}}}

\def\={&\!\!=\!\!&}
\def\bt{\begin{theorem}}
\def\et{\end{theorem}}
\def\bl{\begin{lemma}}
\def\el{\end{lemma}}
\def\br{\begin{remark}}
\def\er{\end{remark}}
\def\bx{\begin{Examples}}
\def\ex{\end{Examples}}
\def\bd{\begin{definition}}
\def\ed{\end{definition}}
\def\bp{\begin{proposition}}
\def\ep{\end{proposition}}
\def\bc{\begin{corollary}}
\def\ec{\end{corollary}}

\def\geq{\geqslant}
\def\leq{\leqslant}

\def\div{\mathord{{\rm div}}}

\def\bP{{\mathbf P}}

\def\<{\langle} \def\>{\rangle}

\allowdisplaybreaks

\begin{document}

\title{Superposition principle for non-local Fokker-Planck operators}

\author{Michael R\"ockner, Longjie Xie and Xicheng Zhang}

\address{Michael R\"ockner:
Fakult\"at f\"ur Mathematik, Universit\"at Bielefeld,
33615, Bielefeld, Germany\\
Email: roeckner@math.uni-bielefeld.de
 }

\address{Longjie Xie:
	School of Mathematics and Statistics, Jiangsu Normal University,
	Xuzhou, Jiangsu 221000, P.R.China\\
	Email: longjiexie@jsnu.edu.cn
}

\address{Xicheng Zhang:
School of Mathematics and Statistics, Wuhan University,
Wuhan, Hubei 430072, P.R.China\\
Email: XichengZhang@gmail.com
 }

\thanks{
This work is supported by NNSF of China (No. 11731009, 11931004), NSF of Jiangsu (BK20170226)
 and the DFG through the CRC 1283 
``Taming uncertainty and profiting from randomness and low regularity in analysis, stochastics and their applications''.
}

\begin{abstract}
	
We prove the superposition principle for probability measure-valued solutions to non-local Fokker-Planck equations, 
which in turn yields the equivalence between martingale problems for SDEs with jumps and such non-local PDEs  with rough coefficients. 
As an application, we obtain a  probabilistic representation for weak solutions of fractional porous media equations.
	
	\bigskip
	
	
\end{abstract}

\maketitle

\tableofcontents

\section{Introduction}

\subsection{Background}

Let $\cP(\mR^d)$ be the space of all probability measures on $\mR^d$ endowed with the weak convergence topology.
Let $b:\mR_+\times\mR^d\to\mR^d$ be a measurable vector field. 
In \cite{Am}, Ambrosio studied the connection between the continuity equation
\begin{align}\label{ce}
\p_t\mu_t=\div(b\mu_t),
\end{align}
and the ordinary differential equation (ODE for short)
\begin{align}\label{ode}
\dif \omega_t=b_t(\omega_t)\dif t.
\end{align}
The following superposition principle was proved therein: Suppose that $t\mapsto\mu_t\in\cP(\mR^d)$ is a solution of (\ref{ce}) and satisfies
$$
\int^T_0\!\!\int_{\mR^d}\frac{|b_t(x)|}{1+|x|}\mu_t(\dif x)\dif t<\infty,\ \ \forall T>0,
$$
then there exists a probability measure $\eta$ on the space $\mC$ of continuous functions from $\mR_+$ to $\mR^d$, 
which is concentrated on the set of all $\omega$ such that $\omega$ is an absolutely continuous solution of (\ref{ode}),
and for every function $f \in C_b(\mR^d)$ and all $t\geq 0$, 
$$
\int_{\mR^d} f(x)\mu_t(\dif x)=\int_{\mC} f\big(\omega_t\big)\eta(\dif \omega).
$$
In other words, the measure $\mu_t$ coincides with the image of $\eta$ under the evaluation map  $\omega\mapsto\omega_t$. 
Consequently, the well-posedness of ODE (\ref{ode}) is equivalent to the existence and uniqueness of solutions for the continuity equation (\ref{ce}).
In particular, the well-posedness of ODE (\ref{ode}) with BV drifts whose distributional divergence belongs to $L^\infty$ was obtained in a generalized sense.
See also \cite{Am2,AFFGP, AT,ST} and the references therein for further developments.

\vspace{2mm}
The stochastic counterpart of the above superposition principle was established by Figalli \cite{Fi}. 
In this situation, the  continuity equation becomes the Fokker-Planck equation, while the ODE becomes
 a stochastic differential equation (SDE for short). 
 More precisely, let $X_t$ solve the following  SDE in $\mR^d$:
\begin{align}\label{sdeb}
\dif X_t=b_t(X_t)\dif t+\sigma_t(X_t)\dif W_t,
\end{align}
where  $b: \mR_+\times\mR^{d}\to\mR^{d}$  and $\sigma: \mR_+\times\mR^{d}\to\mR^{d}\otimes\mR^{d}$ are measurable functions, $W_t$ is a standard Brownian motion defined on some probability space $(\Omega,\cF,\bP)$.
Let $\mu_t\in\cP(\mR^d)$ be the marginal law of $X_t$. By It\^o's formula,
$\mu_t$ solves the following Fokker-Planck equation  in the distributional sense
\begin{align}\label{lfpk}
\p_t\mu_t=\big(\sA_t+\sB_t\big)^*\mu_t,
\end{align}
where for $f\in C_b^2(\mR^d)$,
\begin{align}\label{ab}
\sA_t f(x):=\tr (a_t(x)\cdot\nabla^2 f(x)),\ \sB_t f(x):=b_t(x)\cdot\nabla f(x)
\end{align}
with $a_t(x)=\frac{1}{2}(\sigma_t\sigma^T_t)(x)$, and  $\sA_t^*$ and $\sB_t^*$ stand for the adjoint operators of $\sA_t $ and $\sB_t$, respectively.
When the coefficients $a$ and $b$ are {\it bounded} measurable, the  superposition principle for equation (\ref{lfpk}) 
was proved by Figalli \cite[Theorem 2.6]{Fi}, which says that  every probability measure-valued solution to the Fokker-Planck equation (\ref{lfpk}) 
yields a martingale solution for the operator $\sA_t+\sB_t$ on the path space $\mC$ 
(or equivalently,  a weak solution for SDE (\ref{sdeb})). 
Later, Trevisan \cite{Tr} extended it to the following natural {\it integrability} assumption:
\begin{align}\label{Int}
\int^T_0\!\!\int_{\mR^d}\Big(|b_t(x)|+|a_t(x)|\Big)\mu_t(\dif x)\dif t<\infty,\ \ \forall T>0.
\end{align}
More precisely, for any probability measure-valued solution $\mu$ of \eqref{lfpk}, under \eqref{Int}, there is a weak solution $X$ to SDE \eqref{sdeb}
so that for each $t>0$,
\begin{align}\label{MT1}
\mu_t=\mbox{Law of $X_t$.}
\end{align}
It should be noticed that if $\mu_t$ does not have finite first order moment, then \eqref{Int} may not be satisfied for $b$ and $\sigma$ with at most linear growth.
Recently,  in \cite{BRS}, Bogachev, R\"ockner and  Shaposhnikov obtained the superposition principle under the following more natural assumption:
\begin{align*}
\int^T_0\!\!\int_{\mR^d}\frac{|\<x,b_t(x)\>|+|a_t(x)|}{1+|x|^2}\mu_t(\dif x)\dif t<\infty,\ \ \forall T>0.
\end{align*}
The proofs in \cite{BRS} depend on quite involved uniqueness results for Fokker-Planck equations obtained in \cite{BKRS}.
The superposition principle obtained in \cite{Fi,Tr} has been used in the study of the uniqueness of FPEs with rough coefficients (see e.g. \cite{Ro-Zh, Zh3}), probabilistic representations for solutions to non-linear partial differential equations (PDEs for short) \cite{BM1} as well as distribution dependent SDEs (see  \cite{BM2,R-Z}).  

\medskip

On the other hand, let $(X_t)_{t\geq 0}$ be a Feller process in $\mR^d$ with infinitesimal generator $(\sL, \text{Dom}(\sL))$ (see \cite[page 88]{Re-Yo}).
One says  that $\sL$ satisfies a positive maximum principle if for all
$0\leq f\in \text{Dom}(\sL)$ reaching a positive maximum at point $x_0\in\mR^d$, then $\sL f(x_0)\leq 0$. 
Suppose that $C^\infty_c(\mR^d)\subset\text{Dom}(\sL)$. The well-known Courr\`ege theorem states that
$\sL$ satisfies the positive maximum principle if and only if $\sL$ takes the following form
\begin{align}\label{e:1.1}
\begin{split}
\sL f(x)&=\sum_{i,j=1}^d a_{ij}(x)\p ^2_{ij} f(x)+\sum_{i=1}^{d}b_i(x)\p_if(x)+c(x)f(x)\\
&+\int_{\mR^d}\left(f(x+z)-f(x)-\1_{|z|\leq 1}z\cdot\nabla f(x)\right)\nu_x(\dif z),  
\end{split}
\end{align}
where $a=(a_{ij})_{1\leq i,j\leq d}$ is a $d\times d$-symmetric positive definite matrix-valued measurable function on $\mR^d$,
$b: \mR^d\to\mR^d$, $c:\mR^d\to(-\infty,0]$ are measurable functions and $\nu_x(\dif z)$
is a family of L\'evy measures (see \cite{Sch}). In particular, if we let $\mu_t$ be the marginal law of $X_t$, then by Dynkin's formula,
$$
\p_t\mu_t=\sL^*\mu_t.
$$
We naturally ask that for any probability measure-valued solution $\mu_t$ to the above Fokker-Planck equation, 
is it possible to find some process $X$ so that $\mu_t$ is just the law of $X_t$ for each $t\geq 0$?
In the next subsection, under some growth assumptions on the coefficients, we shall give an affirmative answer.

\subsection{Superposition principle for non-local operators}
Our aim in this paper is to develp a {\it non-local} version of the superposition principle. 
Let $\{\nu_{t,x}\}_{t\geq 0,x\in\mR^d}$ be a family of L\'evy measures over $\mR^d$, that is, for each $t\geq 0$ and $x\in\mR^d$,
\begin{align}\label{GT}
g^\nu_t(x):=\int_{B_\ell}|z|^2\nu_{t,x}(\dif z)<\infty,\quad \nu_{t,x}(B^c_\ell)<\infty,
\end{align}
where $\ell>0$ is a fixed number, and $B_\ell:=\{z\in\mR^d: |z|<\ell\}$. Without loss of generality we may assume
$$
\ell\leq1/\sqrt{2}.
$$
We introduce the following L\'evy type operator: for any $f\in C^2_b(\mR^d)$,
\begin{align}\label{aa2}
\sN_t f(x):=\sN^\nu_t f(x):=\sN^{\nu_{t,x}} f(x):=\int_{\mR^d}\Theta_f(x;z)\nu_{t,x}(\dif z),
\end{align}
where
\begin{align}\label{CA2}
\Theta_f(x;z):=f(x+z)-f(x)-\1_{|z|\leq \ell}\cdot\nabla f(x).
\end{align}
Let us consider the following non-local Fokker-Planck equation (FPE for short):
\begin{align}\label{FPE}
\p_t\mu_t=\sL_t^*\mu_t,
\end{align}
where $\sL_t$ is a general diffusion operator with jumps, i.e.,
\begin{align*}
\sL_t:=\sA_t+\sB_t+\sN_t
\end{align*}
with $\sA_t$ and $\sB_t$ being defined by (\ref{ab}) and $\sN_t$ being defined by \eqref{aa2}. 
We introduce the following definition of weak solution to equation (\ref{FPE}).

\bd[Weak solution]\label{Def11}
Let $\mu:\mR_+\to\cP(\mR^d)$ be a continuous curve. We call $\mu_t$ a weak solution of 
 the non-local FPE (\ref{FPE})
if for any $R>0$ and $t>0$,
\begin{align}\label{CO1}
\left\{
\begin{aligned}
&\int^t_0\!\!\int_{\mR^d}\1_{B_R}(x)\Big(|a_s(x)|+|b_s(x)|+g^\nu_s(x)\Big)\mu_s(\dif x)\dif s<\infty,\\
&\int^t_0\!\!\int_{\mR^d}\Big(\nu_{s,x}(B^c_{\ell\vee(|x|-R)})+\1_{B_R}(x)\nu_{s,x}(B^c_\ell)\Big)\mu_s(\dif x)\dif s<\infty,
\end{aligned}
\right\}
\end{align}
and for all $f\in C^2_c(\mR^d)$ and $t\geq 0$,
\begin{align}\label{ID}
\mu_t(f)=\mu_0(f)+\int^t_0\mu_s(\sL_s f)\dif s,
\end{align}
where $\mu_t(f):=\int_{\mR^d}f(x)\mu_t(\dif x)$.
\ed

We  point out that  unlike the local case considered in \cite{Am,BRS,Fi,Tr}, where the local integrability of the coefficients
with respect to $\mu_t(\dif x)\dif t$ implies the well-definedness of the integrals in \eqref{ID}, 
it is even not clear whether  the above integral in (\ref{ID}) makes sense in the {\it non-local} 
case since in general $\sN^\nu_tf$ does not have compact support for $f\in C^2_c(\mR^d)$. 
This is the reason why we need the second assumption in \eqref{CO1}.

\br\rm
Under \eqref{CO1}, one has $\int^t_0\mu_s(|\sL_s f|)\dif s<\infty$. Let us only show
$$
\int^t_0\mu_s(|\sN^\nu_s f|)\dif s<\infty.
$$
Note that for $x,z\in\mR^d$, by Taylor's expansion, there is a $\theta\in[0,1]$ such that
\begin{align}\label{Ta}
f(x+z)-f(x)-z\cdot\nabla f(x)=z_iz_j\p_i\p_jf(x+\theta z)/2.
\end{align}
Suppose that the support of $f$ is contained in a ball $B_R$. By definition we have
\begin{align*}
|\Theta_f(x;z)|&\leq \|f\|_\infty\1_{|z|>\ell}(\1_{|x+z|<R}+\1_{|x|<R})+\|\nabla^2f\|_\infty\1_{|z|\leq\ell}|z|^2\1_{|x|<R+\ell}.
\end{align*}
Hence,
\begin{align*}
\int^t_0\mu_s(|\sN^\nu_s f|)\dif s&\lesssim \int^t_0\!\!\!\int_{\mR^d}\Big[\nu_{s,x}(B^c_{\ell\vee(|x|-R)})+\1_{B_{R}}(x)\nu_{s,x}(B^c_\ell)\Big]\mu_s(\dif x)\dif s\\
&\quad+\int^t_0\!\!\!\int_{\mR^d}\1_{B_{R+\ell}}(x)g^\nu_s(x)\mu_s(\dif x)\dif s<\infty.
\end{align*}
\er

Let $\mD$ be the space of all $\mR^d$-valued c\'adl\'ag functions  on $\mR_+$, which is endowed with the Skorokhod topology 
so that $\mD$ becomes a Polish space. Let $X_t(\omega)=\omega_t$ be the canonical process. For $t\geq 0$, let $\cB^0_t(\mD)$ 
denote the natural filtration generated by $(X_s)_{s\in[0,t]}$, and let 
$$
\cB_t:=\cB_t(\mD):=\cap_{s\geq t}\cB^0_t(\mD),\ \ \cB:=\cB(\mD):=\cB_\infty(\mD).
$$
Now we recall the notion of martingale solutions associated with $\sL_t$ in the sense of Stroock-Varadhan \cite{St-Va}. 
\bd[Martingale Problem]\label{Def21}
Let $\mu_0\in\cP(\mR^d)$, $s\geq 0$ and $\tau\geq s$ be a $\cB_t$-stopping time.
We call a probability measure $\mP\in\cP(\mD)$ a  martingale solution (resp. a ``stopped'' martingale solution) 
of $\sL_t$ with initial distribution $\mu_0$ at time $s$ if 
\begin{enumerate}[(i)]
\item $\mP(X_t=X_s, t\in[0,s])=1$ and $\mP\circ X_s^{-1}=\mu_0$.

\item For any $f\in C^2_c(\mR^d)$, $M^f_t$ (resp. $M^f_{t\wedge\tau}$) is a $\cB_t$-martingale under $\mP$, where
\begin{align}\label{MMF}
M^f_t:=f(X_t)-f(X_s)-\int^t_s\sL_r f(X_r)\dif r,\ t\geq s.
\end{align}
\end{enumerate}
All the martingale solutions (resp. ``stopped'' martingale solutions) associated with $\sL_t$ with initial law $\mu_0$ at time $s$ will be denoted by $\cM^{\mu_0}_s(\sL)$
(resp. $\cM^{\mu_0}_{s,\tau}(\sL)$).
In particular, if $\mu_0=\delta_x$ (the Dirac measure concentrated on $x$), we shall write $\cM^{x}_s(\sL)=\cM^{\delta_x}_s(\sL)$ for simplify.
\ed

\br\rm \label{Re14}
Under \eqref{Assu1} below, by suitable localization technique, (ii) in Definition \ref{Def21}
is equivalent to that for any $f\in C^2(\mR^d)$ with $|f(x)|\leq C\log(2+|x|)$, $M^f_t$ is a local $\cB_t$-martingale under $\mP$.
\er
  
Throughout this paper, we make the following assumption:
\begin{align}\label{Assu1}
\Gamma^\nu_{a,b}:=\sup_{t,x}\left[\frac{|a_t(x)|+g^\nu_t(x)}{1+|x|^2}+\frac{|b_t(x)|}{1+|x|}+\hbar^\nu_t(x)\right]<\infty,
\end{align}
where $g^\nu_t(x)$ is defined by \eqref{GT} and
\begin{align}\label{GT2}
\hbar^\nu_t(x):=\int_{B_\ell^c}\log\left(1+\tfrac{|z|}{1+|x|}\right)\nu_{t,x}(\dif z),
\end{align}
and if $\nu_{t,x}$ is symmetric, then we define
\begin{align}\label{sym}
\hbar^\nu_t(x):=\int_{|z|>1+|x|}\log\left(1+\tfrac{|z|}{1+|x|}\right)\nu_{t,x}(\dif z).
\end{align}
The main result of this paper is as follows.

\bt[Superposition principle]\label{Main}
Under \eqref{Assu1}, for any weak solution $(\mu_t)_{t\geq 0}$ of FPE \eqref{FPE} in the sense of Definition \ref{Def11},
there is a martingale solution $\mP\in\cM^{\mu_0}_0(\sL_t)$ such that
$$
\mu_t=\mP\circ X^{-1}_t,\ \ \forall t\geq 0.
$$
\et


\br\label{br1}\rm
Under \eqref{Assu1}, condition \eqref{CO1} holds. In fact, it suffices to check that
\begin{align}\label{DQ1}
\sup_{t,x}\Big(\nu_{t,x}(B^c_{\ell\vee(|x|-R)})+\1_{B_R}(x)\nu_{t,x}(B^c_\ell)\Big)<\infty, \ \forall R>0.
\end{align}
By definition we have
\begin{align*}
\nu_{t,x}(B^c_{\ell\vee(|x|-R)})&\leq \int_{B_\ell^c}\log\left(1+\tfrac{|z|}{1+|x|}\right)/\log\left(1+\tfrac{\ell\vee(|x|-R)}{1+|x|}\right)\nu_{t,x}(\dif z)\\
&=h^\nu_{t}(x)/\log\left(1+\tfrac{\ell\vee(|x|-R)}{1+|x|}\right)\leq h^\nu_{t}(x)/\log\left(1+\tfrac{\ell}{1+\ell+R}\right),
\end{align*}
and
\begin{align*}
\1_{B_R}(x)\nu_{t,x}(B^c_\ell)&\leq \1_{B_R}(x)\int_{B_\ell^c}\log\left(1+\tfrac{|z|}{1+|x|}\right)/\log\left(1+\tfrac{\ell}{1+|x|}\right)\nu_{t,x}(\dif z)\\
&=\1_{B_R}(x) h^\nu_{t}(x)/\log\left(1+\tfrac{\ell}{1+|x|}\right)\leq h^\nu_{t}(x)/\log\left(1+\tfrac{\ell}{1+R}\right).
\end{align*}
Hence, \eqref{DQ1} follows by \eqref{Assu1}.
\er
\bx\label{br2}\rm
Let $\nu_{t,x}(\dif z)=\kappa_t(x,z)\dif z/|z|^{d+\alpha}$ with $\alpha\in(0,2)$, that is, $\sN_t$ is an $\alpha$-stable like operator.

\vspace{1mm}
\noindent (i) If $|\kappa_t(x,z)|\leq c(1+|x|)^{\alpha\wedge 1}/(1+\1_{\alpha=1}\log(1+|x|))$, 
then $\sup_{t,x}\hbar^\nu_t(x)<\infty.$
Indeed, by definition we have
\begin{align*}
\hbar^\nu_t(x)\lesssim
\frac{(1+|x|)^{\alpha\wedge 1}}{1+\1_{\alpha=1}\log(1+|x|)}\int_{B_\ell^c}\log\left(1+\tfrac{|z|}{1+|x|}\right)\frac{\dif z}{|z|^{d+\alpha}}.
\end{align*}
We calculate the right hand integral which is denoted by $\sI$ as follows: using polar coordinates and integration by parts,
\begin{align*}
\sI&=c\int^\infty_\ell\log\left(1+\tfrac{r}{1+|x|}\right)r^{-1-\alpha}\dif r\\
&\lesssim \log\left(1+\tfrac{\ell}{1+|x|}\right)+\int^\infty_\ell r^{-\alpha}\left(1+|x|+r\right)^{-1}\dif r\\
&\lesssim (1+|x|)^{-1}+(1+|x|)^{-1}\int^{1+|x|}_\ell r^{-\alpha}\dif r+\int^\infty_{1+|x|}r^{-1-\alpha}\dif r\\
&\lesssim (1+|x|)^{-1}+(1+|x|)^{-(\alpha\wedge1)}(1+\1_{\alpha=1}\log(1+|x|))+(1+|x|)^{-\alpha}\\
&\lesssim (1+|x|)^{-(\alpha\wedge1)}(1+\1_{\alpha=1}\log(1+|x|)).
\end{align*}
Thus, we have $\hbar^\nu_t(x)\leq C$. 

\vspace{1mm}
\noindent (ii)  If $\kappa_t(x,z)$ is symmetric, that is, $\kappa_t(x,z)=\kappa_t(x,-z)$, and $|\kappa_t(x,z)|\leq c(1+|x|)^{\alpha}$, $\alpha\in(0,2)$.
Then $\sup_{t,x}\hbar^\nu_t(x)<\infty.$ In fact, by (\ref{sym}) we have for any $\beta\in(0,\alpha\wedge1)$,
\begin{align*}
\hbar^\nu_t(x)&\lesssim(1+|x|)^{\alpha}\int_{|z|>1+|x|}\left(1+\tfrac{|z|}{1+|x|}\right)^\beta\frac{\dif z}{|z|^{d+\alpha}}\\
&\lesssim(1+|x|)^{\alpha-\beta}\int_{|z|>1+|x|}\frac{\dif z}{|z|^{d+\alpha-\beta}},
\end{align*}
which in turn yields $\sup_{t,x}\hbar^\nu_t(x)<\infty.$
\ex

As far as we know, there are very few results concerning the superposition principle for non-local operators. 
In the constant non-local case,  the third author of the present paper \cite{Zh3} used the superposition principle to show the uniqueness of non-local FPEs. 
Recently, Fournier and Xu \cite{F-X} proved a non-local version to the superposition principle in a special case, that is, 
$$
\sN^\nu_t f(x)=\int_{\mR^d}[f(x+z)-f(x)]\nu_{t,x}(\dif z),
$$
and $(\mu_t)_{t\geq 0}$ have finite first order moments, i.e.,
$$
\int_{\mR^d}|x|\mu_t(\dif x)<\infty,\ \ \forall t\geq 0.
$$
These two assumptions rule out the interesting $\alpha$-stable processes (see Example \ref{br2} above).
To drop these two limitations, we employ some techniques from \cite{BRS}.
It should be emphasized that the elegant push-forward method used in \cite{Tr} does not seem to work in the non-local case. 
Here the main obstacles are to show the tightness and taking limits.
One important motivation for studying the superposition principle for nonlocal operators is 
to solve the Boltzman equation as explained in Subsection 1.2 of \cite{F-X} (see also \cite{HK}).

\subsection{Equivalence between FPEs and martingale problems}
The following corollary is a direct consequence of Theorem \ref{Main} and \cite[Theorem 4.4.2]{EK} (see also \cite[Lemma 2.12]{Tr}). 
For the readers' convenience, we provide a detailed proof here.

\bc\label{Cor1}
Under (\ref{Assu1}), the well-posedness of the Fokker-Planck equation (\ref{FPE}) is 
equivalent to the well-posedness of the martingale problem associated with $\sL$. More precisely,
we have the following equivalences:
\medskip
\begin{enumerate}[$\bullet$]
\item{\bf (Existence)} For any $\nu\in\cP(\mR^d)$, the non-local FPE \eqref{FPE} admits a solution $(\mu_t)_{t\geq 0}$
with initial value $\mu_0=\nu$ if and only if $\cM^{\nu}_0(\sL)$ has at least one element.
\medskip
\item{\bf (Uniqueness)} The following two statements are equivalent.
\begin{enumerate}[(i)]
\item For each $(s,\nu)\in\mR_+\times\cP(\mR^d)$, the non-local FPE \eqref{FPE} has at most one solution $(\mu_t)_{t\geq s}$ with $\mu_s=\nu$
after time $s$.
\item For each $(s,\nu)\in\mR_+\times\cP(\mR^d)$, $\cM^\nu_s(\sL)$ has at most one element.
\end{enumerate}
\end{enumerate}
\ec
\begin{proof}
We only prove the uniqueness part.  (ii)$\Rightarrow$(i) is easy by Theorem \ref{Main}. We show (i)$\Rightarrow$(ii).
For given $(s,\nu)\in\mR_+\times\cP(\mR^d)$ and let $\mP_1,\mP_2\in\cM^\nu_s(\sL)$. 
To show $\mP_1=\mP_2$, it suffices to prove the following claim by induction: 

{\bf (C$_n$)} for given $n\in\mN$, and for any $s\leq t_1<t_2<t_n$ and strictly positive and bounded measurable functions $f_1,\cdots, f_n$ on $\mR^d$,
\begin{align}\label{GK1}
\mE^{\mP_1}(f_1(X_{t_1})\cdots f_n(X_{t_n}))=\mE^{\mP_2}(f_1(X_{t_1})\cdots f_n(X_{t_n})).
\end{align}
First of all, by Theorem \ref{Main} and the assumption, one sees that {\bf (C$_1$)} holds.
Next we assume {\bf (C$_n$)} holds for some $n\geq 2$. For simplicity we write
$$
\eta:=f_1(X_{t_1})\cdots f_n(X_{t_n}),
$$
and for $i=1,2$, we define new probability measures
$$
\dif\tilde\mP_i:=\eta\dif\mP_i/\int_\Omega\eta\dif\mP_i\in\cP(\mD),\ \ \tilde\nu_i:=\tilde\mP_i\circ X^{-1}_{t_n}\in\cP(\mR^d).
$$
Now we show
$$
\tilde\mP_i\in\cM^{\tilde\nu_i}_{t_n}(\sL),\ \ i=1,2.
$$
Let $M^f_t$ be defined by \eqref{MMF}.
We only need to prove that for any $t'>t\geq t_n$ and bounded $\cB_t$-measurable $\xi$,
$$
\mE^{\tilde\mP_i}\left(M^f_{t'}\xi\right)=\mE^{\tilde\mP_i}\left(M^f_t\xi\right)\Leftrightarrow \mE^{\mP_i}(M^f_{t'}\xi\eta)=\mE^{\mP_i}(M^f_t\xi\eta),
$$
which follows since $\mP_i\in\cM^{\nu}_s(\sL)$. Thus, by induction hypothesis and Theorem \ref{Main},
$$
\tilde\nu_1=\tilde\nu_2\Rightarrow\tilde\mP_1\circ X^{-1}_{t_{n+1}}=\tilde\mP_2\circ X^{-1}_{t_{n+1}},\ \ \forall t_{n+1}>t_n.
$$
which in turn implies that {\bf (C$_{n+1}$)} holds. The proof is complete.
 \end{proof}
 
\subsection{Fractional porous media equation}
Probabilistic representation of solution to PDEs is a  powerful tool to study their analytic
properties (well-posedness, regularity, etc) since it allows us to use many probabilistic tools (see \cite{BM2}, \cite{BRR}, \cite{BR}).
As an application of the superposition principle obtained in Theorem \ref{Main}, we intend to derive a probabilistic representation for the weak solution of the following fractional porous media equation (FPME for short):
\begin{align}\label{pde}
\p_tu=\Delta^{\alpha/2}(|u|^{m-1}u),\quad u(0,x)=\varphi(x),
\end{align}
where the porous media exponent $m>1$, $\alpha\in(0,2)$ and $\Delta^{\alpha/2}:=-(-\Delta)^{\alpha/2}$ is the usual fractional Laplacian
with, up to a constant, alternative expression
\begin{align}\label{Fr}
\Delta^{\alpha/2} f(x)={\rm P.V.}\int_{\mR^d}(f(x+z)-f(x))\dif z/|z|^{d+\alpha},
\end{align}
where P.V. stands for the Cauchy principle value.
This equation is a typical non-linear, degenerate and non-local parabolic equation, which  appears naturally in statistical mechanics and population dynamics in order to describe the hydrodynamic limit of interacting particle systems with jumps or long-range interactions. In the last decade, there are many works devoted to the study of equation (\ref{pde}) from the PDE point of view, see \cite{Re-Ro-Wa} 
and the recent survey paper \cite{V}, the monograph \cite{Va} and the references therein. 

\medskip
Let  $\dot{H}^{\alpha/2}(\mR^d)$ be the homogeneous fractional Sobolev space defined as the completion of $C_0^\infty(\mR^d)$
with respect to
$$
\|f\|_{\dot{H}^{\alpha/2}}:=\left(\int_{\mR^d}|\xi|^\alpha|\hat f(\xi)|^2\dif\xi\right)^{1/2}=\|(-\Delta)^{\alpha/4}f\|_2,
$$
where $\hat f$ is the Fourier transform of $f$.
The following notion about the weak solution of FPME is introduced in \cite[Definition 3.1]{PQRV}.
\bd\label{Def10}
A function $u$ is called a weak or $L^1$-energy solution of FPME (\ref{pde}) if
\begin{itemize}
	\item $u\in C([0,\infty);L^1(\mR^d))$ and $|u|^{m-1}u\in L^2_{loc}((0,\infty);\dot{H}^{\alpha/2}(\mR^d))$;
	\item for every $f\in C_0^1(\mR_+\times\mR^d)$, 
	$$
	\int_0^\infty\!\!\int_{\mR^d}u\cdot\p_tf\dif x\dif t=\int_0^\infty\!\!\int_{\mR^d}(|u|^{m-1}u)\cdot \Delta^{\alpha/2}f\dif x\dif t;
	$$
        \item $u(0,x)=\varphi(x)$ almost everywhere.
\end{itemize}
\ed

The following result was proved in \cite[Theorem 2.1, Theorem 2.2]{PQRV}.

\bt\label{t1}
Let $\alpha\in(0,2)$ and $m>1$. For every $\varphi\in L^1(\mR^d)$, 
there exists a unique weak solution $u$ for equation (\ref{pde}). Moreover, $u$ enjoys the following properties:
\begin{enumerate}[(i)]
	\item  if $\varphi\geq 0$, then $u(t,x)>0$ for all $t>0$ and $x\in\mR^d$;

	\item $\p_tu\in L^\infty((s,\infty);L^1(\mR^d))$ for every $s>0$;

	\item for all $t\geq 0$, $\int_{\mR^d}u(t,x)\dif x=\int_{\mR^d}\varphi(x)\dif x$;
	
	\item  if $\varphi\in L^\infty(\mR^d)$, then for every $t>0$,
	\begin{align*}
	\|u(t,\cdot)\|_\infty\leq \|\varphi\|_\infty;
	\end{align*}
	\item for some $\beta\in(0,1)$, $u\in C^\beta((0,\infty)\times\mR^d)$.
\end{enumerate}
\et

Our aim in this subsection is to represent the above solution $u$  as the distributional density of the solution to 
a nonlinear stochastic differential equation driven by the $\alpha$-stable process $L_t$ with L\'evy measure $\dif z/|z|^{d+\alpha}$. 
More precisely, consider the following distribution dependent stochastic differential equation (DDSDE for short) driven by the
$d$-dimensional isotropic $\alpha$-stable process $L_t$:
\begin{align}\label{ddsde}
\dif Y_t&=\rho_{Y_t}\big(Y_{t-}\big)^{\frac{m-1}{\alpha}}\dif L_t,\ \ \rho_{Y_0}(x)=\varphi(x),
\end{align}
where $\rho_{Y_t}(x):=(\dif\cL_{Y_t}/\dif x)(x)$ denotes the distributional density of $Y_t$ with respect to Lebesgue measure. 
We introduce the following notion about the above DDSDE \eqref{ddsde}.
\bd
Let $(\Omega,\cF,\bP; (\cF_t)_{t\geq 0})$ be a stochastic basis  and $(Y,L)$ two $\cF_t$-adapted c\`adl\`ag processes. 
For $\mu\in\cP(\mR^d)$,
we call $(\Omega,\cF,\bP; (\cF_t)_{t\geq 0}; Y,L)$ a solution of \eqref{ddsde} with initial law $\mu$ if
\begin{enumerate}[(i)]
\item $L$ is an $\alpha$-stable process with L\'evy measure $\dif z/|z|^{d+\alpha}$;
\item for each $t\geq 0$, $\bP\circ Y^{-1}_t(\dif x)=\rho_{Y_t}(x)\dif x$;
\item $Y_t$ solves the following SDE:
$$
Y_t=Y_0+\int^t_0\rho_{Y_s}\big(Y_{s-}\big)^{\frac{m-1}{\alpha}}\dif L_s.
$$
\end{enumerate}
\ed
The following is the  second main result of this paper.

\bt\label{main1}
Let $\varphi\geq 0$ be bounded and satisfy $\int_{\mR^d}\varphi(x)\dif x=1$. Let $u$ be the unique weak solution to FPME (\ref{pde}) given by Theorem \ref{t1}
with initial value $\varphi$. Then there exists a weak solution  $Y$ to DDSDE (\ref{ddsde}) such that
\begin{align*}
\rho_{Y_t}(x)=u(t,x),\quad\forall t\geq 0.
\end{align*}
\et
\br\rm
Here an open question is to show the uniqueness of weak solutions to the nonlinear SDE \eqref{ddsde}, 
which can not be derived from the uniqueness of FPME \eqref{Fr}. We will study this in a future work.
\er
We mention that in the 1-dimensional case, such kind of probabilistic representation for the classical porous media equation (i.e., $\alpha=2$) 
was obtained in \cite{BRR}, see also \cite{BRR2} and \cite{BM1, BM2} and  for the generalization to the multi-dimensional case
and more general non-linear equations.  We also mention that there has been an increasing interest in DDSDEs driven 
by Brownian motion in the last decade, see \cite{BM2,R-Z,Wa} and in particularly, \cite{Ca-De1} as well as the references therein.
As far as we know, even the weak existence result  for DDSDE (\ref{ddsde}) driven by L\'evy noise in Theorem \ref{main1} is also new.

\medskip

This paper is organized as follows: In Section 2, we study the equation (\ref{FPE}) with smooth and non-degenerate coefficients.
Then we prove Theorem \ref{Main} and Theorem \ref{main1} in Sections 3 and 4, respectively. 
Throughout this paper we shall use the following conventions:
\begin{itemize}
	\item The letter $C$ denotes a constant, whose value may change in different places.
	\item We use $A\lesssim B$  to denote $A\leq C B$ for some unimportant constant $C>0$.
	\item $\mN_0:=\mN\cup\{0\}$, $\mR_+:=[0,\infty)$, $a\vee b:=\max(a,b)$, $a\wedge b:=\min(a,b)$, $a^+:=a\vee 0$.
	\item $\nabla_x:=\p_x:=(\p_{x_1},\cdots,\p_{x_d})$, $\p_i:=\p_{x_i}:=\p/\p x_i$.
	\item $\mS^d_+$ is the set of all $d\times d$-symmetric and non-negative definite matrices.
\end{itemize}

\section{Proof of Theorem \ref{Main}: Smooth and nondegenerate coefficients}

First of all, we show the following well-posedness result about the martingale problem associated with $\sL_t$,
which extends Stroock's result \cite{St} to unbounded coefficients case, and is probably well-known at least to experts. However, since we can not
find it in the literature, we provide a detailed proof here.
\bt\label{Th21}
Suppose that the following conditions are satisfied:
\begin{enumerate}[{\bf (A)}]
\item $a_t(x):\mR_+\times\mR^d\to\mS^d_+$ is continuous and $a_t(x)$ is invertible;
\item $b_t(x):\mR_+\times\mR^d\to\mR^d$ is locally bounded and measurable;
\item for any $A\in\cB(\mR^d)$, $(t,x)\mapsto\int_A(1\wedge|z|^2)\nu_{t,x}(\dif z)$ is continuous;
\item the following global growth condition holds:
$$
\bar\Gamma^\nu_{a,b}:=\sup_{t,x}\left(\frac{|a_t(x)|+\<x,b_t(x)\>^++g^\nu_t(x)}{1+|x|^2}+2\hbar^\nu_t(x)\right)<\infty,
$$
where $g^\nu_t(x)$ and $\hbar^\nu_t(x)$ are defined by \eqref{GT} and \eqref{GT2}, respectively.
\end{enumerate}
Then for each $(s,x)\in\mR_+\times\mR^d$, there is a unique martingale solution $\mP_{s,x}\in\cM^{x}_s(\sL_t)$.
Moreover, the following assertions hold:
\begin{enumerate}[(i)]
\item For each $A\in\cB(\mD)$, $(s,x)\mapsto \mP_{s,x}(A)$ is Borel measurable. 
\item The following strong Markov property holds: for every $f\in C_b(\mR_+\times\mR^{d})$ and any finite stopping time $\tau$,
$$
\mE^{\mP_{0,x}} (f(\tau+t, X_{\tau+t})|\cB_\tau)=\mE^{\mP_{\tau,X_\tau}} (f(s+t, X_{s+t})).
$$
\end{enumerate}
\et

\br\rm
Condition {\bf (D)} ensures the non-explosion of the solution. 
\er

To prove this theorem we first show the following Lyapunov's type estimate.
\bl
Let $\psi\in C^2(\mR;\mR_+)$ with $\lim_{r\to\infty}\psi(r)=\infty$ and 
\begin{align}\label{CA1}
0< \psi'\leq 1,\quad\psi''\leq 0.
\end{align}
Fix $y\in\mR^d$ and define a Lyapunov function $V_{y}(x):=\psi(\log(1+|x-y|^2))$.
Then for all $t\geq 0$ and $x\in\mR^d$, we have
\begin{align}\label{ag}
\sL_tV_{y}(x)\leq 2\left(\frac{|a_t(x)|+\<x-y,b_t(x)\>^++g^\nu_t(x)}{1+|x-y|^2}+2H^\nu_t(x,y)\right),
\end{align}
where $g^\nu_t(x)$ is defined by \eqref{GT}, and
\begin{align}\label{HH1}
H^\nu_t(x,y):=\int_{B^c_\ell}\log\left(1+\tfrac{|z|}{1+|x-y|}\right)\nu_{t,x}(\dif z).
\end{align}
\el
\begin{proof}
By definition, it is easy to see that
\begin{align*}
\nabla V_y(x)=\frac{2(x-y)}{1+|x-y|^2}\psi'(\log(1+|x-y|^2))
\end{align*}
and
\begin{align*}
\begin{split}
\nabla^2 V_y(x)&=\frac{4(x-y)\otimes (x-y)}{(1+|x-y|^2)^2}(\psi''-\psi')(\log(1+|x-y|^2))\\
&\quad+\frac{2\mI}{1+|x-y|^2}\psi'(\log(1+|x-y|^2)).
\end{split}
\end{align*}
Thus by \eqref{CA1}, one gets that
\begin{align*}
	\sA^a_t V_y(x)\leq \frac{2|a_t(x)|}{1+|x-y|^2},\quad \sB^b_t V_y(x)\leq \frac{2\<x-y,b_t(x)\>^+}{1+|x-y|^2}.
\end{align*}
On the other hand, recalling \eqref{CA2}, we have for $|z|\leq\ell\leq1/\sqrt{2}$,
\begin{align*}
\Theta_{V_y}(x;z)&=V_y(x+z)-V_y(x)-z\cdot \nabla V_y(x)=z_iz_j\p_i\p_jV_y(x+\theta z)/2\\
&=\frac{2\<z,x-y+\theta z\>^2}{(1+|x-y+\theta z|^2)^2}(\psi''-\psi')(\log(1+|x-y+\theta z|^2))\\
&\quad+\frac{|z|^2}{1+|x-y+\theta z|^2}\psi'(\log(1+|x-y+\theta z|^2))\\
&\stackrel{\eqref{CA1}}{\leq}\frac{|z|^2}{1+|x-y+\theta z|^2}\leq\frac{|z|^2}{1+|x-y|^2/2-|z|^2}\leq \frac{2|z|^2}{1+|x-y|^2},
\end{align*}
where $\theta\in[0,1]$. Similarly, by the mean value formula, we have
\begin{align*}
V_y(x+z)-V_y(x)&=\psi'(\theta_*)\Big[\log\big(1+|x-y+z|^2\big)-\log\big(1+|x-y|^2\big)\Big]\\
&\leq\log\left(1+\frac{2|\<x-y,z\>|+|z|^2}{1+|x-y|^2}\right)\leq\log\left(1+\frac{|z|}{\sqrt{1+|x-y|^2}}\right)^2\\
&\leq\log\left(1+\frac{2|z|}{1+|x-y|}\right)^2\leq\log\left(1+\frac{|z|}{1+|x-y|}\right)^4,
\end{align*}
where $\theta_*\in\mR$. Hence,
\begin{align*}
&\sN^\nu_tV_y(x)\leq\int_{\mR^d}\Theta_{V_y}(x;z)\nu_{t,x}(\dif z)\leq 2\frac{g^\nu_t(x)}{1+|x-y|^2}+4H^\nu_t(x,y).
\end{align*}
Combining the above calculations, we obtain \eqref{ag}.
\end{proof}

The following stochastic Gronwall inequality for continuous martingales was proved by Scheutzow \cite{Sc}, and
for general discontinuous martingales in \cite[Lemma 3.7]{XZ2}.

\bl[Stochastic Gronwall inequality]\label{Gron}
Let $\xi(t)$ and $\eta(t)$ be two non-negative c\`adl\`ag adapted processes, 
$A_t$ a continuous non-decreasing adapted process with $A_0=0$, $M_t$ a  local martingale with $M_0=0$. Suppose that
$$
\xi(t)\leq\eta(t)+\int^t_0\xi(s)\dif A_s+M_t,\ \forall t\geq 0.
$$
Then for any $0<q<p<1$ and stopping time $\tau>0$, we have
\begin{align*}
\big[\mE(\xi(\tau)^*)^{q}\big]^{1/q}\leq \Big(\tfrac{p}{p-q}\Big)^{1/q}\Big(\mE \e^{pA_{\tau}/(1-p)}\Big)^{(1-p)/p}\mE\big(\eta(\tau)^*\big),  \label{gron}
\end{align*}
where $\xi(t)^*:=\sup_{s\in[0,t]}\xi(s)$.
\el

The following localization lemma is well known (see e.g. \cite[Theorem 1.3.5]{St-Va}). 
Although it is only proved for the probability measures on the space of continuous functions, by checking the proof therein, one sees that it also works for $\mD$.

\bl\label{Ext}
Let $(\mP_n)_{n\in\mN}\subset\cP(\mD)$ be a family of probability measures and
$(\tau_n)_{n\in\mN}$ a non-decreasing sequence of stopping times with $\tau_0\equiv 0$. 
Suppose that for each $n\in\mN$, $\mP_{n}$ equals $\mP_{n-1}$ on $\cB_{\tau_{n-1}}(\mD)$, 
and for any $T\geq 0$, 
\begin{align*}
\lim_{n\to\infty} \mP_n (\tau_n\leq T)=0.
\end{align*}
Then there is a unique probability measure $\mP\in\cP(\mD)$ such that  $\mP$ equals $\mP_n$ on $\cB_{\tau_n}(\mD)$ and
$\mP_n$ weakly converges to $\mP$ as $n\to\infty$.
\el

We now use the above localization lemma to give
\begin{proof}[Proof of Theorem \ref{Th21}]
Let $\chi\in C^\infty_c(\mR^d)$ be a smooth funciton with 
$$
\chi(x)=1,\ \ |x|<1,\ \ \chi(x)=0,\ \ |x|>2.
$$ 
For any $n\in\mN$, define
$$
\chi_n(x):=\chi(x/n)
$$
and
$$
a^n_t(x):=a_t(x\chi_n(x)),\ b^n_t(x):=\chi_n(x)b_t(x),\ \nu^n_{t,x}(\dif z):=\chi_n(x)\nu_{t,x}(\dif z).
$$
By the assumptions {\bf (A)}-{\bf (C)}, one can check that $(a^n,b^n,\nu^n)$ satisfies for any $T>0$,
\begin{enumerate}[{\bf (A$'$)}]
\item $a^n_t(x):[0,T]\times\mR^d\to\mS^d_+$ is bounded continuous and $a^n_t(x)$ is invertible.
\item $b^n_t(x):[0,T]\times\mR^d\to\mR^d$ is bounded measurable.
\item For any $A\in\cB(\mR^d)$, $(t,x)\mapsto\int_A(1\wedge|z|^2)\nu^n_{t,x}(\dif z)$ is bounded continuous.
\end{enumerate}
Let $\sL^n_t$ be defined in terms of $(a^n,b^n,\nu^n)$. For each $n\in\mN$ and $(s,x)\in\mR_+\times\mR^d$,
by \cite[Theorem 2.34, p.159]{J-S}, there is a unique martingale solution $\mP^n_{s,x}\in\cM^x_s(\sL^n_t)$,
 and the following properties hold:
\begin{enumerate}[(i)]
\item For each $A\in\cB(\mD)$, $(s,x)\mapsto \mP^n_{s,x}(A)$ is Borel measurable. 
\item The following strong Markov property holds: for any $f\in C_b(\mR^{d+1})$ and finite stopping time $\tau$,
$$
\mE^{\mP^n_{0,x}} (f(\tau+t, X_{\tau+t})|\cB_\tau)=\mE^{\mP^n_{\tau,X_\tau}} (f(s+t, X_{s+t})).
$$
\end{enumerate}
Moreover, if we define
$$
\tau_n:=\inf\{t\geq s: |X_t|>n\},
$$
then by \cite[Theorem 2.41, p.161]{J-S}, for any $m\geq n$,
the ``stopped'' martingale problem $\cM^x_{s,\tau_n}(\sL^{m}_t)$ admits a unique solution, that is,
$$
\mP^{m}_{s,x}|_{\cB_{\tau_n}(\mD)}=\mP^{n}_{s,x}|_{\cB_{\tau_n}(\mD)}.
$$
To show the well-posedness, by Lemma \ref{Ext},
it suffices to show that for any $T>0$,
$$
\lim_{n\to\infty}\mP^n_{s,x}(\tau_n\leq T)=0.
$$
Let $V(x):=\log(1+|x|^2)$. By the definition of martingale solution (see Remark \ref{Re14}), there is a c\'adl\'ag local $\mP^n_{s,x}$-martingale $M_t$ such that
\begin{align*}
V(X_{t\wedge\tau_n})&=V(x)+\int^{t\wedge\tau_n}_s\sL^n_rV(X_r)\dif r+M_t\\
&=V(x)+\int^{t\wedge\tau_n}_s\sL_rV(X_r)\dif r+M_t\\
&\stackrel{\eqref{ag}}{\leq} V(x)+4\,\Gamma^{\nu}_{a,b}\cdot(t-s)+M_t.
\end{align*}
where $\Gamma^\nu_{a,b}$ is defined by (\ref{Assu1}).
By Lemma \ref{Gron} and condition {\bf (D)}, we obtain
$$
\sup_n\mE^{\mP^n_{s,x}}\left(\sup_{t\in[s,T\wedge\tau_n]}V^{\frac{1}{2}}(X_t)\right)<+\infty,
$$
which in turn implies that
\begin{align*}
\mP^n_{s,x}(\tau_n\leq T)&=\mP^n_{s,x}\left(\sup_{t\in[s,T\wedge\tau_n]}|X_t|>n\right)\\
&\leq\frac{1}{V^{\frac{1}{2}}(n)}\mE^{\mP^n_{s,x}}\left(\sup_{t\in[s,T\wedge\tau_n]}V^{\frac{1}{2}}(X_t)\right)\to 0 
\end{align*}
as $n\to\infty$. The proof is complete.
\end{proof}

Now we can give the proof of Theorem \ref{Main} under the assumptions {\bf (A)}-{\bf (D)}.
\bt\label{Th24}
Assume that {\bf (A)}-{\bf (D)} hold. Then for any $\mu_0\in\cP(\mR^d)$, there are a unique solution $(\mu_t)_{t\geq 0}$ to FPE \eqref{FPE} and
a unique martingale solution
$\mP_{0,\mu_0}\in\cM^{\mu_0}_0(\sL)$ so that $\mu_t=\mP_{0,\mu_0}\circ X^{-1}_t$.
\et
\begin{proof}
Let $\mu_0\in\cP(\mR^d)$ and $\mP_{s,x}\in\cM^x_s(\sL)$.  Clearly, 
$$
\mP_{0,\mu_0}:=\int_{\mR^d}\mP_{0,x}\mu_0(\dif x)\in\cM^{\mu_0}_0(\sL),
$$
and $\mu_t:=\mP_{0,\mu_0}\circ X^{-1}_t$ solves FPE \eqref{FPE}. It remains to show the uniqueness for  \eqref{FPE}.
Following the same argument as in \cite{F-X}, due to Horowitz and Karandikar \cite[Theorem B1]{HK}, we only need to verify the following five points:
\begin{enumerate}[{\bf (a)}]
\item $C^2_c(\mR^d)$ is dense in $C_0(\mR^d)$ with respect to the uniform convergence.
\item $(t,x)\to\sL_t f(x)$ is measurable for all $f\in C^2_c(\mR^d)$.
\item For each $t\geq 0$, the operator $\sL_t$ satisfies the maximum principle.
\item There exists a countable family $(f_k)_{k\in\mN}\subset C^2_c(\mR^d)$ such that for all $t\geq 0$,
$$
\{\sL_tf, f\in C^2_c(\mR^d)\}\subset\overline{\{\sL_tf_k, k\in\mN\}},
$$
where the closure is taken in the uniform norm.
\item For each $x\in\mR^d$, $\cM^x_0(\sL)$ has exactly one element.
\end{enumerate}
Note that {\bf (a)}-{\bf (c)} are obvious and {\bf (e)} is proven in Theorem \ref{Th21}.
Thus we only need to check {\bf (d)}. Let $(f_k)_{k\in\mN}$ be a countable dense subset of $C^2_c(\mR^d)$, that is, for any $f\in C^2_c(\mR^d)$
with support in $B_R$, where $R\geq 2$, there is a subsequence $f_{k_n}$ with support in $B_{2R}$ such that
$$
\lim_{n\to\infty}\Big(\|f_{k_n}-f\|_\infty+\|\nabla f_{k_n}-\nabla f\|_\infty+\|\nabla^2f_{k_n}-\nabla^2f\|_\infty\Big)=0.
$$
We want to show
$$
\lim_{n\to\infty}\|\sL_t (f_{k_n}-f)\|_\infty=0.
$$
Without loss of generality, we may assume $f=0$ and proceed to prove the following limits:
$$
\lim_{n\to\infty}\|\sA_t f_{k_n}\|_\infty=0,\ \lim_{n\to\infty}\|\sB_t f_{k_n}\|_\infty=0,\ \lim_{n\to\infty}\|\sN^\nu_t f_{k_n}\|_\infty=0.
$$
The first two limits are obvious. Let us focus on the last one. By definition we have
\begin{align*}
|\Theta_{f_{k_n}}(x;z)|&=|f_{k_n}(x+z)-f_{k_n}(x)-\1_{|z|\leq\ell} z\cdot\nabla f_{k_n}(x)|\\
&\leq \1_{|z|>\ell}|f_{k_n}(x+z)|+\1_{|z|>\ell}\1_{B_{2R}}(x)\|f_{k_n}\|_\infty\\
&\quad+\1_{|z|\leq \ell}\1_{B_{2R+2\ell}}(x)\|\nabla^2 f_{k_n}\|_\infty|z|^2.
\end{align*}
Note that
$$
\1_{|z|>\ell}\1_{B_{5R}}(x)\leq \Big[\log(1+\tfrac{\ell}{1+5R})\Big]^{-1}\log(1+\tfrac{|z|}{1+|x|}),
$$
and if $|x|>5R$, then for $|x+z|\leq 2R$, 
$$
\tfrac{|z|}{1+|x|}\geq\tfrac{|x|-|x+z|}{1+|x|}\geq \tfrac{|x|-2R}{1+|x|}>\tfrac{1}{2},
$$ 
and thus,
$$
\1_{|z|>\ell}\1_{B^c_{5R}}(x)\1_{B^c_{2R}}(x+z)\leq \Big[\log(\tfrac{3}{2})\Big]^{-1}\log(1+\tfrac{|z|}{1+|x|}).
$$
Therefore, 
\begin{align*}
|\sN^\nu_tf_{k_n}(x)|&\leq\int_{\mR^d}|\Theta_{f_{k_n}}(x;z)|\nu_{t,x}(\dif z)
\leq \|\nabla^2 f_{k_n}\|_\infty\sup_{x\in B_{2R+2}}\int_{B_1}|z|^2\nu_{t,x}(\dif z)\\
&\quad+C\|f_{k_n}\|_\infty\int_{B^c_1}\log(1+\tfrac{|z|}{1+|x|})\nu_{t,x}(\dif z)\\
&=\|\nabla^2 f_{k_n}\|_\infty\sup_{x\in B_{2R+2}} g^\nu_t(x)+C\|f_{k_n}\|_\infty\sup_{x\in\mR^d}\hbar^\nu_t(x),
\end{align*}
which in turn implies by \eqref{Assu1} that
$$
\lim_{n\to\infty}\|\sN^\nu_tf_{k_n}\|_\infty=0.
$$
The proof is compete.
\end{proof}

\section{Proof of Theorem \ref{Main}: General case}

Let $\mu_t$ be a solution of \eqref{FPE} in the sense of Definition \ref{Def11}. In order to show the existence of a martingale solution $\mP\in\cM^{\mu_0}_0(\sL_t)$ so that
$$
\mu_t=\mP\circ X^{-1}_t,
$$
we shall follow the same lines of argument as in \cite{Fi} and \cite{Tr}. Here and below we use the following convention: for $t\leq 0$,
$$
\mu_t(\dif x):=\mu_0(\dif x),\ a_t(x)=0,\ b_t(x)=0,\ \nu_{t,x}(\dif z)=0.
$$
\subsection{Regularization}
Let $\rho^{\rm t}\in C_c^\infty([0,1];\mR_+)$ with $\int^1_0\rho^{\rm t}(s)\dif s=1$ and 
$\rho^{\rm x}\in C_c^\infty(B_1;\mR_+)$ with $\int_{\mR^d}\rho^{\rm x}(x)\dif x=1$. For $\eps>0$, define
$$
\rho^{\rm t}_\eps(t):=\eps^{-1}\rho(t/\eps),\ \rho^{\rm x}_\eps(x):=\eps^{-d}\rho(x/\eps),\ \rho_\eps(t,x):=\rho^{\rm t}_\eps(t)\rho^{\rm x}_\eps(x).
$$
Given a locally finite signed measure $\zeta_t(\dif x)\dif t$ on $\mR^{d+1}$, we define
$$
\rho_\eps*\zeta(t,x):=\int_{\mR^{d+1}}\rho_\eps(t-s,x-y)\zeta_s(\dif y)\dif s.
$$
Throughout this section we shall fix 
$$
\ell\in(0,1/\sqrt{2}).
$$
We first show the following regularization estimate.
\bl\label{Le41}
Let $a,b$ and $\nu$ be as in the introduction. For  $\eps\in(0,\ell)$, we have
\begin{align*}
\frac{|\rho_\eps*(a\mu)|(t,x)}{1+|x|^2}&\leq \sup_{t,y}\frac{2|a_t(y)|}{1+|y|^2}(\rho_\eps*\mu)(t,x),\\
\frac{|\rho_\eps*(b\mu)|(t,x)}{1+|x|}&\leq \sup_{t,y}\frac{2|b_t(y)|}{1+|y|}(\rho_\eps*\mu)(t,x).
\end{align*}
Moreover, if we let 
\begin{align*}
\bar{\nu}^\eps_{t,x}(\dif z):=\int_{\mR^{d+1}}\rho_\eps(t-s,x-y)\nu_{s,y}(\dif z)\mu_s(\dif y)\dif s,
\end{align*}
then we also have 
\begin{align*}
\frac{g^{\bar{\nu}^\eps}_t(x)}{1+|x|^2}&\leq \sup_{t,y}\frac{2g^\nu_{t}(y)}{1+|y|^2}(\rho_\eps*\mu)(t,x),\\
H^{\bar{\nu}^\eps}_t(x,y)&\leq 2\sup_{t,x}H^{\nu}_t(x,y)(\rho_\eps*\mu)(t,x),
\end{align*}
where $g^\nu_{t}(x)$ and $H^{\nu}_t(x,y)$ are defined by \eqref{GT} and \eqref{GT2}, respectively.
\el
\begin{proof}
Note that for $|x-y|\leq \ell\leq1/\sqrt{2}$,
\begin{equation}\label{DT2}
(1+|y|^2)/2\leq 1+|x|^2\leq 2(1+|y|^2).
\end{equation}
Fix $\eps\in(0,\ell)$ below. By definition we have
\begin{align*}
\frac{|\rho_\eps*(a\mu)|(t,x)}{1+|x|^2}&\leq \int_{\mR^{d+1}}\rho_\eps(t-s, x-y)\frac{|a_s(y)|}{1+|x|^2}\mu_s(\dif y)\dif s\\
&\leq 2\int_{\mR^{d+1}}\rho_\eps(t-s,x-y)\frac{|a_s(y)|}{1+|y|^2}\mu_s(\dif y)\dif s,
\end{align*}
and
\begin{align*}
\frac{|\rho_\eps*(b\mu)|(t,x)}{1+|x|}&\leq\int_{\mR^{d+1}}\rho_\eps(t-s,x-y)\frac{|b_s(y)|}{1+|x|}\mu_s(\dif y)\dif s\\
&\leq 2\int_{\mR^{d+1}}\rho_\eps(t-s,x-y)\frac{|b_s(y)|}{1+|y|}\mu_s(\dif y)\dif s.
\end{align*}
Similarly, by Fubini's theorem and \eqref{DT2}, we have
\begin{align*}
\frac{g^{\bar{\nu}^\eps}_t(x)}{1+|x|^2}&=\int_{\mR^{d+1}}\int_{B_\ell}\frac{|z|^2}{1+|x|^2}\rho_\eps(t-s,x-y)\nu_{s,y}(\dif z)\mu_s(\dif y)\dif s\\
&\leq 2\int_{\mR^{d+1}}\int_{B_\ell}\frac{|z|^2}{1+|y|^2}\rho_\eps(t-s, x-y)\nu_{s,y}(\dif z)\mu_s(\dif y)\dif s\\
&=2\int_{\mR^{d+1}}\frac{g^\nu_s(y)}{1+|y|^2}\rho_\eps(t-s,x-y)\mu_s(\dif y)\dif s,
\end{align*}
and
\begin{align*}
H^{\bar{\nu}^\eps}_t(x,y)&=\int_{\mR^{d+1}}\int_{B_\ell^c}\log\left(1+\tfrac{|z|}{1+|x-y|}\right)\rho_\eps(t-s,x-y')\nu_{s,y'}(\dif z)\mu_s(\dif y')\dif s\\
&\leq\int_{\mR^{d+1}}\int_{B_\ell^c}\log\left(1+\tfrac{2|z|}{1+|y'-y|}\right)\rho_\eps(t-s,x-y')\nu_{s,y'}(\dif z)\mu_s(\dif y')\dif s\\
&\leq 2\int_{\mR^{d+1}}H^{\nu}_s(y',y)\rho_\eps(t-s,x-y')\mu_s(\dif y')\dif s.
\end{align*}
Combining the above calculations, we obtain the desired estimates.
\end{proof}
Let  $\phi(x):=(2\pi)^{-d}\e^{-|x|^2/2}$ be the normal distribution density. For $\eps\in(0,\ell)$, 
as in \cite{BRS}, we define the approximation sequence $\mu^\eps_t\in\cP(\mR^d)$ by
\begin{align}\label{mo}
\mu_t^\eps(x):=(1-\eps)(\rho_\eps*\mu)(t,x)+\eps\phi(x).
\end{align}
We have the following easy consequence. 
\bp
\begin{enumerate}[(i)]
\item For each $t\geq 0$ and $\eps\in(0,\ell)$, we have
$$
0<\mu^\eps_t(x)\in C^\infty(\mR_+; C^\infty_b(\mR^{d})),\ \ \int_{\mR^d}\mu_t^\eps(x)\dif x=1.
$$
\item For each $t\geq 0$, $\mu^\eps_t$ weakly converges to $\mu_t$, that is, for any $f\in C_b(\mR^d)$,
$$
\lim_{\eps\to\infty}\int_{\mR^d}f(x)\mu^\eps_t(x)\dif x=\int_{\mR^d}f(x)\mu_t(\dif x).
$$
\item $\mu^\eps_t$ solves the following Fokker-Planck equation:
$$
\p_t\mu^\eps_t=(\sA_t^{\eps}+\sB^{\eps}_t+\sN^{\eps}_t)^*\mu^\eps_t=:(\sL_t^\eps)^*\mu^\eps_t,
$$
where $\sA^{\eps}_t$, $\sB^{\eps}_t$ and $\sN^{\eps}_t$ are defined as in the introduction in terms of
\begin{align*}
a^\eps_t(x)&:=\frac{(1-\eps)[\rho_\eps*(a\mu)](t,x)+\eps\phi(x)\mI}{\mu_t^\eps(x)},\\
b^\eps_t(x)&:=\frac{(1-\eps)[\rho_\eps*(b\mu)](t,x)+\eps\phi(x)x}{\mu_t^\eps(x)},
\end{align*}
and
\begin{align}\label{neps}
\nu^\eps_{t,x}(\dif z):=\frac{1-\eps}{\mu_t^\eps(x)}\int_{\mR^{d+1}}\rho_\eps(t-s,x-y)\nu_{s,y}(\dif z)\mu_s(\dif y)\dif s.
\end{align}
\item The following uniform estimates hold: for any $\eps\in(0,\ell)$,
\begin{align}\label{UN0}
\begin{split}
&\sup_{t,x}\left[\frac{|a^\eps_t(x)|+g^{\nu^\eps}_t(x)}{1+|x|^2}
+\frac{|b^\eps_t(x)|}{1+|x|}\right]\leq 1+2\sup_{t,x}\left[\frac{|a_t(x)|+g^{\nu}_t(x)}{1+|x|^2}
+\frac{|b_t(x)|}{1+|x|}\right]
\end{split}
\end{align}
and
\begin{align}\label{UN00}
\sup_{t,x}H^{\nu^\eps}_t(x,y)\leq \sup_{t,x}H^{\nu}_t(x,y),\ \ y\in\mR^d.
\end{align}
\end{enumerate}
\ep
\begin{proof}
The first two assertions are obvious by definition. Let us show (iii). By definition, it suffices to prove that for any 
$f\in C^\infty_c(\mR^d)$ and $t\geq 0$,
\begin{align}\label{DT33}
\mu^\eps_t(f)=\mu^\eps_0(f)+\int^t_0\mu^\eps_s(\sL_s^\eps f)\dif s,
\end{align}
where
$$
\mu^\eps_t(f):=\int_{\mR^d}f(x)\mu^\eps_t(x)\dif x.
$$
Note that for any $f\in C^\infty_c(\mR^d)$,
$$
\Delta\phi+\div(x\cdot \phi)\equiv 0\Rightarrow\int_{\mR^d}\phi(x)(\Delta f(x)-x\cdot\nabla f(x))\dif x=0.
$$
By Fubini's theorem and a change of variables, it is easy to see that \eqref{DT33} holds.
Finally, estimate \eqref{UN0} follows by Lemma \ref{Le41}.
\end{proof}

The following result follows by Theorem \ref{Th24}.

\bl
For any $\eps\in(0,\ell)$ and $(s,x)\in\mR_+\times\mR^d$, there is a unique martingale solution $\mP^\eps_{s,x}\in\cM^x_s(\sL^\eps_t)$.
In particular, there is also a martingale solution $\mQ^\eps\in\cM^{\mu^\eps_0}_0(\sL^\eps_t)$ so that for each $t\geq 0$,
$$
\mu^\eps_t(x)\dif x=\mQ^\eps\circ X^{-1}_t(\dif x).
$$
\el
\begin{proof}
By Theorem \ref{Th24}, it suffices to check that $(a^\eps,b^\eps,\nu^\eps)$ satisfies conditions {\bf (A)}-{\bf (D)}. First of all, {\bf (A)} and {\bf (B)} are obvious,
and {\bf (D)} follows by \eqref{UN0}. It remains to check {\bf (C)}.
We only check that for any $\eps\in(0,\ell)$, $n\in\mN$ and $x,x'\in B_n$, $t,t'\in[0,n]$,
\begin{align}\label{DT3}
\int_{\mR^d}(1\wedge |z|^2)|\nu^\eps_{t,x}-\nu^\eps_{t',x'}|(\dif z)\leq c_{n,\eps}(|t-t'|+|x-x'|).
\end{align}
Noting that 
$$
\inf_t\inf_{x\in B_n}\mu^\eps_t(x)\geq\eps\inf_{x\in B_n}\phi(x),
$$
we have by definition that for all $x,x'\in B_n$ and $t,t'\in[0,n]$,
\begin{align*}
|\nu^\eps_{t,x}-\nu^\eps_{t',x'}|(\dif z)&\leq
\int_{\mR^{d+1}}\left|\frac{\rho_\eps(t-s,x-y)}{\mu_t^\eps(x)}-\frac{\rho_\eps(t'-s,x'-y)}{\mu_{t'}^\eps(x')}\right|\nu_{s,y}(\dif z)\mu_s(\dif y)\dif s\\
&\leq c_{n,\eps}(|t-t'|+|x-x'|)\int_0^{n+1}\int_{B_{n+1}}\nu_{s,y}(\dif z)\mu_s(\dif y)\dif s.
\end{align*}
Estimate \eqref{DT3} then follows since $\sup_{s,y\in[0,n+1]\times B_{n+1}}\int_{\mR^{d}}(1\wedge|z|^2)\nu_{s,y}(\dif z)<\infty$.
\end{proof}

\subsection{Tightness}

We first prepare the following result (cf. \cite[Proposition 7.1.8]{BKRS}).

\bl\label{Le34}
For $\mu^\eps_0\in\cP(\mR^d)$ being defined by \eqref{mo}, there exits a function $\psi\in C^2(\mR_+)$ 
with the properties
\begin{align*}
\psi\geq 0,\quad\psi(0)=0,\quad 0< \psi'\leq 1,\quad -2\leq \psi''\leq 0,\quad\lim_{r\to\infty}\psi(r)=+\infty,
\end{align*}
and such that
\begin{align}\label{JP1}
\sup_{\eps\in[0,\ell)}\int_{\mR^d}\psi\big(\log(1+|x|^2)\big)\mu^\eps_0(\dif x)<\infty.
\end{align}
\el
\begin{proof}
Since $\mu^\eps_0$ weakly converges to $\mu_0$ as $\eps\to0$, we have
$$
\lim_{n\to\infty}\sup_{\eps\in[0,\ell)}\mu^\eps_0(B^c_n)=0.
$$
In particular, we can find a subsequence $n_k$ such that 
for $z_k:=\log(1+n^2_k)$, 
$$
z_{k+1}-z_k\geq z_k-z_{k-1}\geq 1,
$$
and
$$
\sup_{\eps\in[0,\ell)}\int_{\mR^d}1_{[z_k,\infty)}(\log(1+|x|^2))\mu^\eps_0(\dif x)=\sup_{\eps\in[0,\ell)}\mu^\eps_0(B^c_{n_k})\leq 2^{-k}.
$$
Let $z_0=0$ and define 
$$
\psi_0(s):=\sum_{k=0}^\infty\1_{[z_k,z_{k+1}]}(s)\left[k-1+\frac{s-z_k}{z_{k+1}-z_k}\right].
$$ 
Clearly, we have
\begin{align*}
\int_{\mR^d}\psi_0(\log(1+|x|^2))\mu^\eps_0(\dif x)\leq \sum_{k=0}^\infty k\int_{\mR^d}1_{[z_{k},\infty)}(\log(1+|x|^2))\mu^\eps_0(\dif x)
\leq\sum_{k=0}^\infty \frac{k}{ 2^{k}}.
\end{align*}
However, $\psi_0$ does not belong to the class $C^2(\mR_+)$. Let us take 
$$
\psi(t):=\int_0^tg(r)\dif r
$$
with $g\in C^1(\mR_+)$, $0\leq g\leq 1$, $-2\leq g'\leq 0$, and 
$$
g(z)=\psi_0'(z)\quad\text{if}\quad z\in(z_k,z_{k+1}-k^{-1}).
$$
It is easy to see that such a function $g$ always exists. The proof is complete.
\end{proof}

\bl
Let $H^\nu_t(x,y)$ be defined by \eqref{HH1}. We have
\begin{align}\label{HH2}
H^\nu_t(x,y)\leq 2(1+|y|)\hbar^\nu_t(x),\ \ \forall  t\geq 0, x,y\in\mR^d.
\end{align}
\el
\begin{proof}
Recall that
$$
H^\nu_t(x,y)=\int_{B^c_\ell}\log\left(1+\tfrac{|z|}{1+|x-y|}\right)\nu_{t,x}(\dif z).
$$
If $|x|\leq 2|y|$, then
\begin{align*}
H^\nu_t(x,y)&\leq\int_{B^c_\ell}\log\left(1+|z|\right)\nu_{t,x}(\dif z)\leq
\int_{B^c_\ell}\log\left(1+\tfrac{(1+2|y|)|z|}{1+|x|}\right)\nu_{t,x}(\dif z)\\
&\leq\int_{B^c_\ell}\log\left(1+\tfrac{|z|}{1+|x|}\right)^{1+2|y|}\nu_{t,x}(\dif z)=(1+2|y|)\hbar^\nu_t(x).
\end{align*}
If $|x|>2|y|$, then $2|x-y|\geq 2|x|-2|y|\geq |x|$ and
$$
H^\nu_t(x,y)\leq\int_{B^c_\ell}\log\left(1+\tfrac{2|z|}{2+|x|}\right)\nu_{t,x}(\dif z)\leq 2\hbar^\nu_t(x).
$$
The proof is complete.
\end{proof}

Now, we prove the following tightness result.

\bl
The family of probability measures $(\mQ^\eps)_{\eps\in(0,\ell)}$ is tight in $\cP(\mD)$.
\el
\begin{proof}
By Aldous' criterion (see \cite{Al} or \cite[p.356]{J-S}), it suffices to check the following two conditions:
\begin{enumerate}[(i)]
\item For any $T>0$, it holds that
$$
\lim_{N\to\infty}\sup_\eps\mQ^\eps\left(\sup_{t\in[0,T]}|X_t|>N\right)=0.
$$
\item For any $T,\delta_0>0$ and stopping time $\tau<T-\delta_0$, it holds that
$$
\lim_{\delta\to 0}\sup_\eps\sup_\tau\mQ^\eps\left(|X_{\tau+\delta}-X_\tau|>\lambda\right)=0,\ \ \forall\lambda>0.
$$ 
\end{enumerate}

\noindent {\bf Verification of (i)}. Let $\psi$ be as in Lemma \ref{Le34} and $V(x):=\psi(\log(1+|x|^2))$.
By the definition of a martingale solution (see Remark \ref{Re14}), there is a c\'adl\'ag local $\mQ^\eps$-martingale $M^\eps_t$ 
and constant $C$ independent of $\eps$ such that for all $t\geq 0$,
\begin{align*}
V(X_{t})&=V(X_0)+\int^{t}_0\sL^\eps_rV(X_r)\dif r+M_t^\eps\stackrel{\eqref{UN0}}{\leq} V(X_0)+Ct+M^\eps_t.
\end{align*}
By Lemma \ref{Gron} , there is a constant $C>0$ such that for all $T>0$,
\begin{align}\label{EE}
\sup_{\eps\in(0,\ell)}\mE^{\mQ^\eps}\left(\sup_{t\in[0,T]}V^{\frac{1}{2}}(X_t)\right)\leq 
C\sup_{\eps\in(0,\ell)}(\mE^{\mQ_\eps} V(X_0))^{\frac{1}{2}}\stackrel{\eqref{JP1}}{<}\infty,
\end{align}
which in turn implies that (i) is true.

\noindent {\bf Verification of (ii)}. Let $\tau\leq T-\delta_0$ be a bounded stopping time. For any $\delta\in(0,\delta_0)$,
by the strong Markov property  we have
\begin{align}\label{AQ1}
\mQ^\eps\left(|X_{\tau+\delta}-X_\tau|>\lambda\right)&=\mE^{\mQ^\eps}\left(\mP^\eps_{s,y}\left(|X_{s+\delta}-y|>\lambda\right)|_{(s,y)=(\tau, X_\tau)}\right).
\end{align}
Recalling that $V_{y}(x):=\psi(\log(1+|x-y|^2)$, and by \eqref{ag}, \eqref{UN0}, \eqref{UN00}, \eqref{HH2}  and  \eqref{Assu1} we deduce that
\begin{align*}
\sL^\eps_tV_{y}(x)&\leq 2\left(\frac{|a^\eps_t(x)|+\<x-y,b^\eps_t(x)\>^++g^{\nu^\eps}_t(x)}{1+|x-y|^2}+2H^{\nu^\eps}_t(x,y)\right)\\
&\leq C\left(\frac{1+|x|^2+|x-y|(1+|x|)}{1+|x-y|^2}+H^\nu_t(x,y)\right)\leq C(1+|y|^2),
\end{align*}
where $C>0$ is independent of $t,x,y$ and $\eps$. Furthermore, we have
\begin{align*}
V_{y}(X_t)&=V_y(X_0)+\int^{t}_s\sL^\eps_rV_{y}(X_r)\dif r+M^\eps_t\\
&\leq V_y(X_0)+C(1+|y|^2)(t-s)+M^\eps_t,
\end{align*}
where $(M^\eps_t)_{t\geq s}$ is a local $\mP^\eps_{s,y}$-martingale.
By Lemma \ref{Gron} again and since $V_{y}(y)=0$, we obtain
$$
\mE^{\mP^\eps_{s,y}} V_{y}(X_t)^{1/2}\leq C(1+|y|)(t-s)^{1/2}.
$$
Hence,
\begin{align*}
\mP^\eps_{s,y}\left(|X_{s+\delta}-y|>\lambda\right)&=\mP^\eps_{s,y}\left(V_{y}(X_{s+\delta})>\psi(\log(1+\lambda^2))\right)\\
&\leq\mE^{\mP^\eps_{s,y}}\left(V_{y}(X_{s+\delta})\right)/\psi(\log(1+\lambda^2))\\
&\leq C(1+|y|)\delta^{1/2}/\psi(\log(1+\lambda^2)),
\end{align*}
and by \eqref{AQ1} and \eqref{EE},
\begin{align*}
\mQ^\eps\left(|X_{\tau+\delta}-X_\tau|>\lambda\right)&\leq \mQ^\eps(|X_\tau|>R)+C(1+R)\delta^{1/2}/\psi(\log(1+\lambda^2))\\
&\leq C/\psi(\log(1+R^2))+C_\lambda(1+R)\delta^{1/2}.
\end{align*}
Letting $\delta\to 0$ first and then $R\to\infty$, one sees that (ii) is satisfied. 
\end{proof}
\subsection{Limits}
In order to take weak limits, we rewrite 
\begin{align*}
\sB_tf(x)+\sN_tf(x)&=\tilde b_t(x)\cdot\nabla f(x)
+\int_{\mR^d}\Theta^\pi_f(x;z)\nu_{t,x}(\dif z)=:\widetilde{\sB_t}f(x)+\widetilde\sN_tf(x),
\end{align*}
where 
$$
\tilde b_t(x):=b_t(x)+\int_{\mR^d}\big[\pi(z)-z\1_{|z|\leq\ell}\big]\nu_{t,x}(\dif z),
$$
and
\begin{align}\label{npi}
\Theta^\pi_f(x;z):=f(x+z)-f(x)-\pi(z)\cdot\nabla f(x).
\end{align}
Here,
$\pi:\mR^d\to\mR^d$ is a smooth symmetric function satisfying
$$
\pi(z)=z,\ \ |z|\leq\ell,\ \ \pi(z)=0,\ \ |z|>2\ell.
$$
As in (\ref{aa2}), we shall also write $\widetilde\sN_tf(x)=\widetilde\sN_t^\nu f(x)=\widetilde\sN^{\nu_{t,x}}f(x)$.
We have the following result.

\bl\label{Le37}
For any $f\in C^2_c(\mR^d)$ with support in $B_R$, 
there is a constant $C=C(f)>0$ such that for all $x\in\mR^d$ and $z,z'\in \mR^d$ with $|z'|\leq|z|$,
$$
|\Theta^\pi_f(x;z)-\Theta^\pi_f(x;z')|\leq C(|z-z'|\wedge 1)(\1_{B_{R+\ell}}(x)\1_{|z|\leq \ell}|z|+\1_{|z|>\ell\vee(|x|-R)}).
$$
\el
\begin{proof}
Note that
$$
\sQ:=|\Theta^\pi_f(x;z)-\Theta^\pi_f(x;z')|=|f(x+z)-f(x+z')-(\pi(z)-\pi(z'))\cdot\nabla f(x)|.
$$
We make the following decomposition:
\begin{align*}
\sQ=\sQ\cdot\1_{|z|\leq\ell}+\sQ\cdot\1_{|z|>\ell}\1_{|x|\leq R}+\sQ\cdot\1_{|z|>\ell}\1_{|x|>R}=:\sQ_1+\sQ_2+\sQ_3.
\end{align*}
For $\sQ_1$, since supp$(f)\subset B_R$, we have by \eqref{Ta} that
\begin{align*}
|\sQ_1|&\leq|z-z'|^2\|\nabla^2 f\|_\infty\1_{B_{R+\ell}}(x)\1_{|z|\leq\ell}\leq C(|z-z'|\wedge 1)|z|\1_{B_{R+\ell}}(x)\1_{|z|\leq\ell}.
\end{align*}
For $\sQ_2$, we have
\begin{align*}
|\sQ_2|&\leq\Big(|f(x+z)-f(x+z')|+|\pi(z)-\pi(z')|\cdot\|\nabla f\|_\infty\Big)\1_{|z|>\ell}\1_{|x|\leq R}\\
&\leq C(|z-z'|\wedge 1)\1_{|z|>\ell}\1_{|x|\leq R}.
\end{align*}
As for $\sQ_3$, we have
\begin{align*}
|\sQ_3|=|f(x+z)-f(x+z')|\cdot\1_{|z|>\ell}\1_{|x|>R}\leq C(|z-z'|\wedge 1)\1_{|z|>\ell\vee(|x|-R)},
\end{align*}
where we have used that for $|z'|\leq|z|\leq|x|-R$,
$$
f(x+z)=f(x+z')=0.
$$
Combining the above calculations, we obtain the desired estimate.
\end{proof}

 The following approximation result will be crucial for taking weak limits.
\bl\label{imp}
For any $\delta\in(0,1)$ and $T>0$, there is 
a family of L\'evy measures $\eta_{t,x}(\dif z)$ such that for any $f\in C^2_c(\mR^d)$, 
\begin{align}\label{LP2}
\int^T_0\!\!\int_{\mR^d}\sup_{x\in B_1(y)}|\widetilde\sN^{\nu_{s,y}}f(x)-\widetilde\sN^{\eta_{s,y}}f(x)|\mu_s(\dif y)\dif s\leq\delta,
\end{align}
and
$$
\sup_{s,y}\|\widetilde\sN^{\eta_{s,y}}f\|_\infty<\infty,\ \ (s,y,x)\mapsto\widetilde\sN^{\eta_{s,y}}f(x)\mbox{ is continuous}.
$$
\el
\begin{proof}
(i) By the randomization of kernel functions (see \cite[Lemma 14.50, p.469]{Ja}), there is a measurable function 
$$
h_{t,x}(\theta):[0,T]\times\mR^d\times[0,\infty)\to\mR^d\cup\{\infty\}
$$ 
such that
$$
\nu_{t,x}(A)=\int^\infty_0\1_A(h_{t,x}(\theta))\dif\theta,\ \ \forall A\in\sB(\mR^d).
$$
In particular, we have
\begin{align}\label{HQ1}
\widetilde\sN^{\nu_{s,y}} f(x)=\int^\infty_0\Theta^\pi_f(x; h_{s,y}(\theta))\dif\theta=:\widetilde\sN^{h_{s,y}}f(x),
\end{align}
and
$$
g^\nu_t(x)=\int^\infty_0\1_{B_\ell}(h_{t,x}(\theta))|h_{t,x}(\theta)|^2\dif\theta,\ \ 
\nu_{t,x}(B^c_\ell)=\int^\infty_0\1_{B^c_\ell}(h_{t,x}(\theta))\dif\theta.
$$
We introduce $\mX:=[0,T]\times\mR^d\times[0,\infty)$ and a locally finte measure $\gamma$ over $\mX$ by
$$
\gamma(\dif \theta,\dif x,\dif t):= \big(1+\nu_{t,x}(B^c_{\ell\vee(|x|-R)})\big)\dif\theta\mu_{t}(\dif x)\dif t.
$$
{\it Claim:} There is a sequence of measurable functions $\{\bar h^n_{t,x}(\theta), n\in\mN\}$
so that for each $\theta\geq 0$ and $n\in\mN$, $(t,x)\mapsto\bar h^n_{t,x}(\theta)$ is continuous with compact support, and
\begin{align}\label{HQ3}
|\bar h^n_{t,x}(\theta)|\leq |h_{t,x}(\theta)|,
\end{align}
and
\begin{align}\label{HQ4}
\lim_{n\to\infty}\int_\mX\Big(|\bar h^n_{t,x}(\theta)-h_{t,x}(\theta)|^2\wedge 1\Big)\gamma(\dif\theta,\dif x,\dif t)=0.
\end{align}
{\it Proof of Claim:}
Fix $m\in\mN$. Since $\1_{[0,m]}(\theta)\gamma(\dif\theta,\dif x,\dif t)$ is a finite measure (see \eqref{DQ1}) over $\mX$, by Lusin's theorem, there exists  a family 
of continuous function $\{\bar h^\eps_{t,x}(\theta),\eps\in(0,1)\}$ with compact support in $(t,x)$ such that
$$
|\bar h^\eps_{t,x}(\theta)|\leq |h_{t,x}(\theta)|,\ \bar h^\eps_{t,x}(\theta)\to h_{t,x}(\theta),\ \eps\to 0, \gamma-a.s.
$$
Thus by the dominated convergence theorem,
$$
\lim_{\eps\to 0}\int_\mX\Big(|\bar h^\eps_{t,x}(\theta)-h_{t,x}(\theta)|^2\wedge 1\Big)\1_{[0,m]}(\theta)\gamma(\dif\theta,\dif x,\dif t)=0.
$$
On the other hand, we have
$$
\lim_{m\to \infty}\int_\mX\Big(|h_{t,x}(\theta)|^2\wedge 1\Big)\1_{(m,\infty)}(\theta)\gamma(\dif\theta,\dif x,\dif t)=0.
$$
By a diagonalizaion argument, we obtain the desired approximation sequence. The claim is proven.
\medskip
\\
(ii) Let $f\in C^2_c(B_R)$.
By \eqref{HQ1}, \eqref{HQ3} and Lemma \ref{Le37}, we have for all $x\in B_1(y)$,
\begin{align*}
&|\widetilde\sN^{h_{s,y}}f(x)-\widetilde\sN^{\bar h^n_{s,y}}f(x)|\leq\int^\infty_0|\Theta^\pi_f(x; h_{s,y}(\theta))-\Theta^\pi_f(x; \bar h^n_{s,y}(\theta))|\dif\theta\\
&\lesssim\int^\infty_0\Big(|h_{s,y}(\theta)|\1_{B_\ell}(h_{s,y}(\theta))\1_{B_{R+\ell}}(x)+\1_{B^c_{\ell\vee(|x|-R)}}(h_{s,y}(\theta))\Big)\\
&\qquad\qquad\times\Big(|h_{s,y}(\theta)-\bar h^n_{s,y}(\theta)|\wedge 1\Big)\dif\theta\\
&\leq\left(\int^\infty_0\Big(|h_{s,y}(\theta)|^2\1_{B_\ell}(h_{s,y}(\theta))\1_{B_{R+\ell+1}}(y)
+\1_{B^c_{\ell\vee(|x|-R)}}(h_{s,y}(\theta))\Big)\dif\theta\right)^{\frac{1}{2}}\\
&\quad\qquad\times\left(\int^\infty_0\Big(|h_{s,y}(\theta)-\bar h^n_{s,y}(\theta)|^2\wedge 1\Big)\dif\theta\right)^{\frac{1}{2}}\\
&\leq\left(\1_{B_{R+\ell+1}}(y)g^\nu_s(y)+\nu_{s,y}(B^c_{\ell\vee(|y|-R-1)})\right)^{\frac{1}{2}}
\!\left(\int^\infty_0\!\!\Big(|h_{s,y}(\theta)-\bar h^n_{s,y}(\theta)|^2\wedge 1\Big)\dif\theta\right)^{\frac{1}{2}}.
\end{align*}
Hence, by \eqref{Assu1} and \eqref{DQ1} we further have
\begin{align*}
&\int^T_0\!\!\int_{\mR^d}\sup_{x\in B_1(y)}|\sN^{h_{s,y}}f(x)-\sN^{\bar h^n_{s,y}}f(x)|\mu_s(\dif y)\dif s\\
&\lesssim\int^T_0\!\!\int_{\mR^d}\left(\int^\infty_0\Big(|h_{s,y}(\theta)-\bar h^n_{s,y}(\theta)|^2\wedge 1\Big)\dif\theta\right)^{\frac{1}{2}}
\mu_s(\dif y)\dif s\\
&\leq\sqrt{T}\left(\int_\mX(|h_{s,y}(\theta)-\bar h^n_{s,y}(\theta)|^2\wedge 1)\dif\theta\mu_s(\dif y)\dif s\right)^{\frac{1}{2}}\stackrel{\eqref{HQ4}}{\to} 0.
\end{align*}
(iii) For fixed $n\in\mN$, by 
the above claim, it is easy to see that
$$
(s,y,x)\mapsto\sN^{\bar h^n_{s,y}}f(x)\mbox{ is continuous.}
$$
Moreover, we have
\begin{align*}
|\widetilde\sN^{\bar h^n_{s,y}}f(x)|&\leq\int^\infty_0|\Theta^\pi_f(x; \bar h^n_{s,y}(\theta))|\dif\theta
\leq C\int^\infty_0\Big(|\bar h^n_{s,y}(\theta)|^2\wedge 1\Big)\dif\theta.
\end{align*}
Since $\bar h^n_{s,y}(\theta)$ has compact support in $(s,y)$, we have
$$
\sup_{s,y}\|\sN^{\bar h^n_{s,y}}f\|_\infty<\infty.
$$
Finally we just need to take $n$ large enough and define
$$
\eta_{t,x}(A):=\int^\infty_0\1_A(\bar h^n_{t,x}(\theta))\dif\theta.
$$
The proof is complete.
\end{proof}

Now we are in a position to give:

\begin{proof}[Proof of Theorem \ref{Main}]
Let $\mQ$ be any accumulation point of $(\mQ^\eps)$. By taking weak limits for
$$
\mu^\eps_t=\mQ^\eps\circ X^{-1}_t,
$$
we obtain
$$
\mu_t=\mQ\circ X^{-1}_t.
$$
It remains to show that $\mQ\in\cM^{\mu_0}_0(\sL_t)$. We need to show that for any $f\in C^2_c(\mR^d)$,
$$
M_t:=f(X_t)-f(X_0)-\int^t_0\sL_s f(X_s)\dif s
$$
is a $\cB_t$-martingale under $\mQ$. More precisely, it suffices to prove that for any $s<t$ and bounded $\cB_s$-measurable $g_s\in C_b(\mD)$,
$$
\mE^\mQ (M_t g_s)=\mE^\mQ (M_s g_s).
$$
Since $\mQ^\eps\in\sM^{\mu^\eps_0}_0(\sL^\eps_t)$, by the definition of martingale solutions, we have
$$
\mE^{\mQ^\eps} (M^\eps_t g_s)=\mE^{\mQ^\eps} (M^\eps_s g_s),
$$
where 
$$
M^\eps_t:=f(X_t)-f(X_0)-\int^t_0\sL^\eps_s f(X_s)\dif s.
$$
Clearly, we only need to show the following limits.
\begin{align*}
\lim_{\eps\to 0}\mE^{\mQ^\eps}\left (g_s\int^t_s\sA^{\eps}_r f(X_r)\dif r \right)&=\mE^{\mQ} \left(g_s\int^t_s\sA_r f(X_r)\dif r\right),\\
\lim_{\eps\to 0}\mE^{\mQ^\eps}\left (g_s\int^t_s\widetilde{\sB^{\eps}_r} f(X_r)\dif r \right)&=\mE^{\mQ} \left(g_s\int^t_s\widetilde{\sB_r} f(X_r)\dif r\right),\\
\lim_{\eps\to 0}\mE^{\mQ^\eps}\left (g_s\int^t_s\widetilde\sN^{\nu^\eps}_r f(X_r)\dif r \right)&=\mE^{\mQ} \left(g_s\int^t_s\widetilde\sN^\nu_r f(X_r)\dif r\right).
\end{align*}
The first two limits can be proved by following the arguments as in \cite{Fi, Tr}. 
We show the last one, which is more difficult.
Let $\eta_{t,x}(\dif z)$ be as given by  Lemma \ref{imp}, and recall that $\nu^\eps$ is defined by (\ref{neps}). We write
\begin{align*}
&\left|\mE^{\mQ^\eps}\left (g_s\int^t_s\widetilde\sN^{\nu^\eps}_r f(X_r)\dif r \right)-\mE^{\mQ} 
\left(g_s\int^t_s\widetilde\sN^\nu_r f(X_r)\dif r\right)\right|\\
&\leq\left|\mE^{\mQ^\eps}\left (g_s\int^t_s\widetilde\sN^{\nu^\eps}_r f(X_r)\dif r \right)-\mE^{\mQ^\eps} 
\left(g_s\int^t_s\widetilde\sN^{\eta_\eps}_r f(X_r)\dif r\right)\right|\\
&+\left|\mE^{\mQ^\eps}\left (g_s\int^t_s\widetilde\sN^{\eta^\eps}_r f(X_r)\dif r \right)-\mE^{\mQ^\eps} \left(g_s\int^t_s\widetilde\sN^{\eta}_r f(X_r)\dif r\right)\right|\\
&+\left|\mE^{\mQ^\eps}\left (g_s\int^t_s\widetilde\sN^{\eta}_r f(X_r)\dif r \right)-\mE^{\mQ} \left(g_s\int^t_s\widetilde\sN^\eta_r f(X_r)\dif r\right)\right|\\
&+\left|\mE^{\mQ}\left (g_s\int^t_s\widetilde\sN^{\eta}_r f(X_r)\dif r \right)-\mE^{\mQ} \left(g_s\int^t_s\widetilde\sN^\nu_r f(X_r)\dif r\right)\right|
=:\sum_{i=1}^4I_i(\eps),
\end{align*}
where $\eta^\eps$ is defined similarly as in (\ref{neps}) with $\nu$ being replaced by $\eta$.
For $I_1(\eps)$, by  definition, we have
\begin{align*}
I_1(\eps)&\leq\|g_s\|_\infty\mE^{\mQ^\eps}\left(\int^t_s|\widetilde\sN^{\nu^\eps}_r f(X_r)-\widetilde\sN^{\eta^\eps}_r f(X_r)|\dif r\right)\\
&=\|g_s\|_\infty\int^t_s\int_{\mR^d}|\widetilde\sN^{\nu^\eps}_r f(x)-\widetilde\sN^{\eta^\eps}_r f(x)|\mu^\eps_r(x)\dif x\dif r\\
&=(1-\eps)\|g_s\|_\infty\int^t_s\int_{\mR^d}|\widetilde\sN^{\bar\nu^\eps}_r f(x)-\widetilde\sN^{\bar\eta^\eps}_r f(x)|\dif x\dif r\\
&=(1-\eps)\|g_s\|_\infty\int^t_s\int_{\mR^d}\left|\int_{\mR^d}\Theta^\pi_f(x;z)(\bar\nu^\eps_{r,x}-\bar\eta^\eps_{r,x})(\dif z)\right|\dif x\dif r,
\end{align*}
where $\Theta^\pi_f(x;z)$ is defined by \eqref{npi} and 
$$
\bar\nu^\eps_{r,x}(\dif z):=\int_{\mR^{d+1}}\rho_\eps(r-s,x-y)\nu_{s,y}(\dif z)\mu_s(\dif y)\dif s.
$$
By Fubini's theorem we further have
\begin{align*}
I_1(\eps)&\leq\|g_s\|_\infty\int^t_s\int_{\mR^d}\bigg|\int_{\mR^{d+1}}\rho_\eps(r-s,x-y)\\
&\qquad\qquad\qquad\qquad\quad\times\Big(\widetilde\sN^{\nu_{s,y}}f(x)-\widetilde\sN^{\eta_{s,y}}f(x)\Big)\mu_s(\dif y)\dif s\bigg|\dif x\dif r\\
&\leq\|g_s\|_\infty\int^T_0\!\!\!\int_{\mR^d}\sup_{x\in B_1(y)}|\widetilde\sN^{\nu_{s,y}}f(x)-\widetilde\sN^{\eta_{s,y}}f(x)|\mu_s(\dif y)\dif s
\stackrel{\eqref{LP2}}{\leq}\|g_s\|_\infty\delta.
\end{align*}
For $I_2(\eps)$, recalling \eqref{mo},  we have
\begin{align*}
I_2(\eps)&\leq\|g_s\|_\infty\mE^{\mQ^\eps}\left(\int^t_s|\widetilde\sN^{\eta^\eps}_r f(X_r)-\widetilde\sN^{\eta}_r f(X_r)|\dif r\right)\\
&=\|g_s\|_\infty\int^t_s\int_{\mR^d}|\widetilde\sN^{\eta^\eps}_r f(x)-\widetilde\sN^{\eta}_r f(x)|\mu^\eps_r(x)\dif x\dif r\\
&=\|g_s\|_\infty\int^t_s\int_{\mR^d}|(1-\eps)\widetilde\sN^{\bar\eta^\eps}_r f(x)-\mu^\eps_r(x)\widetilde\sN^{\eta}_r f(x)|\dif x\dif r\\
&\leq(1-\eps)\|g_s\|_\infty\int^t_s\int_{\mR^d}\int_{\mR^{d+1}}\rho_\eps(r-s,x-y)\\
&\quad\times|\widetilde\sN^{\eta_{s,y}}f(x)-\widetilde\sN^{\eta_{r,x}}f(x)|\mu_s(\dif y)\dif s\dif x\dif r\\
&\quad+\eps\|g_s\|_\infty\int^t_s\int_{\mR^d}|\phi(x)\widetilde\sN^{\eta}_r f(x)|\dif x\dif r.
\end{align*}
Since $(s,y,x)\mapsto \widetilde\sN^{\eta_{s,y}} f(x)$ is continuous and
$\|\widetilde\sN^\eta f\|_\infty<\infty$, by the dominated convergence theorem, we  get
$$
\lim_{\eps\to 0}I_2(\eps)=0.
$$
Concerning $I_3(\eps)$, it follows by the definition of weak convergence that
$$
\lim_{\eps\to 0}I_3(\eps)=0.
$$
For $I_4(\eps)$, we have
$$
I_4(\eps)\leq \|g_s\|_\infty\int^t_s\int_{\mR^d}|\widetilde\sN^{\eta}_r f(x)-\widetilde\sN^{\nu}_r f(x)|\mu_r(\dif x)\dif r\stackrel{\eqref{LP2}}{\leq}\|g_s\|_\infty\delta.
$$
The proof is complete.
\end{proof}

\section{Proof of Theorem \ref{main1}}

Let   $u$ be the unique weak solution of FPME (\ref{pde}) given by Theorem \ref{t1}
with initial value $\varphi\geq 0$ being bounded and $\int_{\mR^d}\varphi(x)\dif x=1$. Let
$$
\sigma_t(x):=|u(t,x)|^{\frac{m-1}{\alpha}},\ \ \kappa_t(x):=u(t,x)^{m-1},\ \nu_{t,x}(\dif z):=\tfrac{\kappa_t(x)\dif z}{|z|^{d+\alpha}}.
$$
By the change of variable we have
\begin{align}\label{41}
\nu_{t,x}(A)=\int_{\mR^d}\1_{A}(\sigma_t(x) z)\frac{\dif z}{|z|^{d+\alpha}},\ \ A\in\sB(\mR^d\setminus\{0\}),
\end{align}
and
$$
\sN_t f(x):={\rm P.V.}\int_{\mR^d}(f(x+\sigma_t(x)z)-f(x))\frac{\dif z}{|z|^{d+\alpha}}=\kappa_t(x)\Delta^{\alpha/2},
$$
where the second equality is due to \eqref{Fr}.
By Definition \ref{Def10} it is easy to see that $u(t,x)$ solves the following non-local FPE:
$$
\p_t u=\sN^*_t u,\ \ u(0,x)=\varphi(x),
$$
that is, for every $t>0$ and $f\in C_0^2(\mR^d)$,
\begin{align*}
&\int_{\mR^d}f(x)u(t,x)\dif x=\int_{\mR^d}f(x)\varphi(x)\dif x+\int_0^t\!\!\int_{\mR^d}\kappa_s(x) \Delta^{\alpha/2}f(x)u(s,x)\dif x\dif s.
\end{align*}
Note that for each $t>0$, 
$$
|\sigma_t(x)|=|u(t,x)|^{\frac{m-1}{\alpha}}\leq \|\varphi\|_\infty^{\frac{m-1}{\alpha}}.
$$
Thus, by Example \ref{br2} with the above $\nu_{t,x}$ and Theorem \ref{Main} with $\mu_0(\dif x)=\varphi(x)\dif x,$
there is a martingale solution $\mP\in\sM^{\mu_0}_0(\cN_t)$ so that
$$
\mP\circ X^{-1}_t(\dif x)=u(t,x)\dif x,\ \ t\geq 0.
$$ 
By \eqref{41} and \cite[Theorem 2.26, p.157]{J-S}, there are a stochastic basis $(\Omega,\cF,\bP; (\cF_t)_{t\geq 0})$
and a Poisson random measure $N$ on $\mR^d\times  [0,\infty)$ with intensity  $|z|^{-d-\alpha}\dif z\dif t$,
as well as an $\cF_t$-adapted c\`adl\`ag process $Y_t$ such that
$$
\bP\circ Y^{-1}_t(\dif x)=\mP\circ X^{-1}_t(\dif x),\ \ t\geq 0,
$$
and
$$
\dif Y_t=\int_{|z|\leq 1}\sigma_t(Y_{t-})z\tilde N(\dif z,\dif t)+\int_{|z|>1}\sigma_t(Y_{t-})z N(\dif z, \dif t),
$$
where $\tilde N(\dif z,\dif t):=N(\dif z,\dif t)-|z|^{-d-\alpha}\dif z\dif t$.
Finally we just need to define
$$
L_t:=\int_0^t\!\!\int_{|z|\leq 1}z\tilde N(\dif z,\dif s)+\int_0^t\!\!\int_{|z|> 1}z N(\dif z,\dif s),
$$
then $L$ is a $d$-dimensional isotropic 
$\alpha$-stable process with L\'evy measure $\dif z/|z|^{d+\alpha}$, and
$$
\dif Y_t=\sigma_t(Y_{t-})\dif L_t.
$$
The proof is finished.

\end{document}